\documentclass[11pt, a4]{article}

\usepackage{amsmath}
\usepackage{amssymb}
\usepackage{amscd}

\usepackage{url}

\usepackage{endnotes}
\let\footnote=\endnote

\usepackage{graphicx}
\usepackage{float}

\usepackage{array}

\newtheorem{theorem}{Theorem}[section]
\newtheorem{lemma}[theorem]{Lemma}

\newtheorem{proposition}[theorem]{Proposition}

\title{\bf The consistency of arithmetic from a point of view of constructive tableau method with 
strong negation, Part I: the system without complete induction\footnote{This is the second version of the paper, August 13, 2021. One typographical mistake was corrected in the rule 3.6g on page 8. My current affiliations were modified.} }
\author{Takao Inou\'{e}}
\date{}

\begin{document}

\maketitle

\begin{abstract}
In this Part I, we shall prove the consistency of arithmetic without complete induction from a point of view of strong negation, using its embedding to 
the tableau system $\bf SN$ of constructive arithmetic with strong negation without complete induction, for which two types of 
cut elimination theorems hold. 
One is $\bf SN$-cut elimination theorem for the full system $\bf SN$. The other is $\bf PCN$-cut elimination theorem for a proposed subsystem 
$\bf PCN$ of $\bf SN$. The disjunction property and the E-theorem (existence property) for $\bf SN$ are also proved. As a novelty, we 
shall give a simple proof of a restricted version of $\bf SN$-cut elimination theorem as an application of 
the disjunction property, using  $\bf PCN$-cut elimination theorem. 
\end{abstract}

\smallskip



\section{Introduction}

We shall  be concerned in this paper with the consistency of arithmetic without complete induction, namely, Peano's fifth axiom. The consistency of 
arithmetic with or without complete induction has been studied rather intensively by a number of authors such as Ackermann [1924, 1940] 
\cite{ackermann1924, ackermann1940} Neumann [1927] \cite{neumann}, 
Herbrand [1931] \cite{herbrand}, Gentzen [1934 and 1936] \cite{gentzen1934a, gentzen1934b, gentzen1936}, Ono [1938] \cite{onok}, 
and Novikov [1959] \cite{novikov} and in post-war times Hlodovskii (or Kholdovskii) [1959] 
\cite{khlo} and Sch\"utte [1960 and 1977] \cite{schutteB1960, schutte1977}. (Also refer to Hilbert and Bernays [1970, Vol.2, 2nd. ed.] \cite{hb1970}.) 
 (For the details of the recent development, refer to Aczel [1992] \cite{aczel}, Arai [2011] \cite{arai}, 
Buchholz et al [1989] \cite{buch}, Buss [1998] \cite{busshandbook}, Buss and Ignjatovi'c [1995] \cite{bussIg}, Feferman [1981] \cite{feferman}, 
Feferman and Sieg [1981a and 1981b] \cite{fs1, fs2}, Kahle and Rathjen [2015] \cite{kr}, 
Pohlers [2009] \cite{pohlers}, Sch\"utte [1977] \cite{schutte1977}, Takeuti [1975] \cite{takeuti} and Toledo [1975] \cite{toledo}. (For Novikov's proof theory, 
refer to Bellotti \cite{bellotti}. For the history, see Murawski \cite{mura}.) 

All throughout these works, the logical or mathematical system, the consistency of which is at issue, is reduced to (or embedded in) another 
system which is suceptible to the tableau method.As well-known, the system to be developed by the tableau method is cut-free, 
and this property thereof is made use of for proving the consistency.

In what follows, the tableau method to be employed is that for arithmetic based on first-order constructive predicate logic with strong negation. Since the 
proposed predicate logic is constructive, the proposed reduction constitutes a Kolmogorov-Gentzen-G\"odel-type interpretation of classical arithmetic by 
its constructive counterpart (Kolmogorov [1925] \cite{kolmo}, G\"odel[1932-3] \cite{gedel}  and Gentzen [1936] \cite{gentzen1936}). (For tableau methods, 
refer to Fitting \cite{fitting}, Smullyan \cite{smullyan}, etc.)

Now the purpose of this paper is to give a consistency proof of arithmeticc without  complete induction on the basis of the principle as described. 
As already mentioned the constructive system to which classical arithmetic is reduced is not traditional intuitionistic or Heyting arithmetic, but one to be 
developed on the basis of constructive predicate logic with strong negation which takes place of intuitionistic or Heyting negation. As well known, strong 
negation was first introduced by Nelson [1949] \cite{nelson} in connection with recursive realizability. The negation was incorporated into constructive logic 
involving not only strong negation, but also intuitionistic one as well by Markov [1950] \cite{markov}, Vorob'ev [1952 and 1964] \cite{voro1, voro2}, 
Rasiowa [1958] \cite{rasiowa} and others. The proposed arithmetic is based on these studies. The paper to follow consists of ten sections inclusive of 
this introduction. \S 2 and \S 3 consern the preparatories for classical arithmetic and its constructive counterpart $\bf SN$ without complete induction 
to be developed on the basis of constructive predicate logic with strong negation.  (For a recent literature for 
strong negation, refer to Odintsov \cite{odin}.) 

More specifically, classical arithmetic is transformed into a proper subsystem of the cut-free arithmetic, namely, its tableau version, and the consistency of 
the whole system of constructive arithmetic is straightforward as is the case with any cut-free logic. And from this follows the consistency of classical arithmetic. 
As a result, the actual infinity involved, in the intuitive interpretation of classical arithmetic by way of the law of excluded midddle is completely done away with 
by this embedding, since constructive logic does not need any notion of actual infinity.

In \S 4, we shall introduce two subsystems   $\bf PCN$ and  $\bf FN$ of $\bf SN$. Especially, $\bf PCN$ will play essential roles in this paper.

In \S 5, a number of theorems will be proved for facilitating the sebsequent development. The disjunction property and the E-theorem (existence property) 
for $\bf SN$ are also proved. \S 6 will concern two cut elimination theorems, one called $\bf PCN$-cut elimination theorem is indispensable for
 reducing classical arithmetic to its constructive version with the help of the embedding theorem of the former to 
the latter. (Refer to Inou\'e [1984 and 1984] \cite{inouen, inoueac}.) In \S 7, we deal with the reduction, i.e., embedding theorem, which 
constitutes the core of this paper giving rise to the consistency result. In \S 8, we prove the consistency of arithmetic without complete induction. 

In \S 9, we shall present the proof of $\bf SN$-cut elimination theorem for the full system $\bf SN$. We shall more concentrate on the restricted 
version of $\bf SN$-cut elimination theorem, We shall give a simple proof of it, using the disjunction property of $\bf SN$ and 
 $\bf PCN$-cut elimination theorem.

The last \S 10 suggests future studies in the direction of this paper.


\section{Classical arithmetic without complete induction}

Here we wish to present classical arithmetic without complete induction called $\bf CN$ in its Hilbert-type version.

$\bf CN$ could be defined in a number of ways, but here it is defined as an sxiomatic system with the axiom schemata of the following kinds.

2.1. The axiom schemata for classical predicate logic to be developed in terms of six logical symbols, namely, $\wedge$ (conjunction), $\vee$ (disjunction), $\supset$ 
(implication), $\sim$ (classical negation), $\forall$ (universal quantifier), and $\exists$ (existential quantifier):

\medskip 

2.11a $\enspace  \enspace A \supset . B \supset A$,

2.11b $\enspace  \enspace A \supset B . \supset . (A \supset . B \supset C) \supset (A \supset C)$,

2.12 $\enspace  \enspace \enspace A \supset . B \supset (A \wedge B)$,

2.13a $\enspace  \enspace A \wedge B . \supset A$,

2.13b $\enspace  \enspace A \wedge B . \supset B$,

2.14a $\enspace  \enspace A \supset . A \vee B$,

2.14b $\enspace  \enspace B \supset . A \vee B$,

2.15 $\enspace  \enspace A \supset C . \supset . (B \supset C) \supset (A \vee B . \supset C)$, 

2.16 $\enspace  \enspace A \supset B . \supset . (A \supset \sim B) \supset \sim A)$,

2.17 $\enspace  \enspace \sim \sim A \supset A$,

2.18 $\enspace  \enspace \forall x A(x) . \supset A(t)$,

2.19 $\enspace  \enspace A(t) \supset \exists x A(x)$, 

\medskip 

\noindent where $t$ is a term such that no free occurrence of $x$ in $A(x)$ is in the scope of a quantifier $\forall y$ or $\exists y$ with being a variable of $t$.

2.2 The axiom schemata for arithmetic in terms of $'$ (successor function), $;$ (addition), $\cdot$ (multiplication), $=$ (equality):

\medskip

2.21 $\enspace  \enspace a' = b' \supset a = b$, 

2.22 $\enspace  \enspace \sim a' = 0$, 

2.23 $\enspace  \enspace a = b \supset . b = c \supset a = c$, 

2.24 $\enspace  \enspace a = b \supset a' = b'$, 

2.25 $\enspace  \enspace a + 0 = a$, 

2.26 $\enspace  \enspace a + b' = (a + b)' $, 

2.27 $\enspace  \enspace a \cdot 0 = 0$, 

2.28 $\enspace  \enspace a \cdot b' = a \cdot b +a$.

\medskip

Hereby, following Sch\"utte [1977] \cite{schutte1977}, $\bf CN \it \vdash^n A$ ($n \geq 0$) means that $A$ is provable in $\bf CN$ by a proof of the length $n$. 
($A$ is an axiom if $n = 0$.) 

\medskip

2.3 The rules of inference are as follows:

\medskip

2.31 

\bigskip

$\enspace \enspace \displaystyle{\frac{\enspace \mbox{\bf CN} \vdash^{n_1} A \enspace \enspace \enspace 
\mbox{\bf CN} \vdash^{n_2} A \supset B \enspace }{\enspace \enspace  \bf CN \it \vdash^n B \enspace }}$ 

\bigskip

\noindent where $n = max( n_1, n_2) + 1$.

\medskip 

2.32

\bigskip

$\enspace \enspace \displaystyle{\frac{\enspace \mbox{\bf CN} \vdash^n C \supset A(x) \enspace \enspace }
{\enspace \enspace  \mbox{\bf CN} \vdash^{n + 1}  C \supset \forall A(x) \enspace }}$ 

\bigskip

\noindent where $C$ is a formula which does not contain $x$ free.

2.33

\bigskip

$\enspace \enspace \displaystyle{\frac{\enspace \mbox{\bf CN} \vdash^n A(x) \supset C \enspace \enspace }
{\enspace \enspace  \mbox{\bf CN} \vdash^{n + 1} \exists A(x) \supset C \enspace }}$ 

\bigskip

\noindent where $C$ is a formula which does not contain $x$ free. This is an arithmetic taken over from Kleene [1952] \cite{kleene1} if we 
suppress the axiom of complete induction. The well-formed expressions including terms and formulas are defined in the well-known way 
in terms of the primitive symbols. And we do not use different syntactic symbols respectively for free and bound variables. Now, two 
theses of $\bf CN$ will be presented without proof.

\medskip 

2.41 $\enspace \enspace \mbox{\bf CN} \vdash^1 x = x$, 

2.42 $\enspace \enspace \mbox{\bf CN} \vdash^1 x = y \supset . A(x) \supset A(y)$.

\begin{theorem} 
Every general recursive predicate is completely representavle in $\bf CN$.
\end{theorem}

For the proof of Theorem 2.1, refere to Kleene [1952] \cite{kleene1}. 


\section{Constructive arithmetic without complete induction}

Constructive arithmetic without complete induction denoted by $\bf SN$ is now develped by way of the tableau method, or the cut-free  Gentzen-type 
formulation. It is again emphasized that the tableau method is a cut-free system, and it is one of our main results that classical arithmetic in its 
Hilbert- type version is embedded in this cut-free system, of which the consistency is forthcoming outright. Now, we introduce the formal languane 
of $\bf SN$. As primitive symbols we use:

\medskip 

3.11 A denumerably infinite number variables. (Different syntactic symbols are not used respectively for free and bound variables.)

3.12 The symbols $0$ (zero), $'$ (successor function), $\sim$ (strong negation), $\wedge$ (conjunction), $\vee$ (disjunction), $\supset$ 
(implication), $\forall$ (universal quantifier), and $\exists$ (existential quantifier): We note that negation in $\bf SN$ is always regarded as 
strong negation. We shall use the same negation symbol $\sim$ for strong negation as classical one. 

3.13 Symbols for $n$-place general recursive function and $n$-place general recursive predicates ($n \geq 1$). 

3.14 Round brackets and comma.

The terms of $\bf SN$ are defined as usual in terms of $0$, number variables and functionsymbols. $0$ and the expressions to be obtained by 
successively applying successor function $'$ to $0$ are called \it numerals\rm .$0^{(n)}$ stands for 
$$0^{\overbrace{' \cdots '}^{\mbox{\footnotesize \it n \rm times}}},$$
which is the result to be obtained by applying successor function $'$ to $0$ times. A term is \it numerical \rm if it contains no free number 
variables. The prime formulas of $\bf SN$ are $P(t_1, \dots , t_n)$ where $P$ is a symbol for an $n$-place general recursive predicate 
($n \geq 1$) and $t_1, \dots , T_n$ are terms. A prime formula is said to be \it constant \rm if it contains no free number variable. We define the 
\it value \rm of a numerical term as follows.

\medskip
 
3.21 The term $0$  has value $0$.

3.22 If $t$ is a numerical term of value $m$, then $t'$ has value $m'$. 

3.23 If $f$ is a symbol for an $n$-place general recursive function ($n \geq 1$) and $t_1, \dots , t_n$ are terms with values $m_1 , \dots , m_n$, then 
 $f(t_1, \dots , t_n)$ has the value given as $f(m_1, \dots , m_n)$. Clearly every numerical term has a uniquely determined value which is a numeral. 
We now proceed to the definition of the \it truth-value \rm of a constant prime formula in order to introduce the axiom of $\bf SN$.

3.3 If $P$ is a symbol for an $n$-place general recursive predicate ($n \geq 1$) and $t_1 ,\dots , t_n$ are numerical terms with values 
$m_1 , \dots , m_n$, then $P(t_1 , \dots , t_n)$ is \it true \rm or \it false \rm according as $P(m_a , \dots , m_n)$ is decided to be true or false. 
The formula is decidably determined whether it is true or not. 

It is necessary that every constant prime formula is decidably determined whether it is true or not. The formula $A(s_1 , \dots , s_n)$ is said to 
be \it equivalent \rm to $A(r_1, \dots , r_n)$ if $s_1, \dots , s_n$, $r_1, \dots , r_n$ are numerical terms with  $s_1, \dots , s_n$, respectively, 
having the same values as $r_1, \dots , r_n$. It is immediate that equivalence of formulas is an equivalence relation. (Outermost round branckets 
will be suppressed wherever no ambiguity  arises therefrom.) This language is essentially due to Sch\"utte [1977] 
\cite[p.169 for the system $\Delta^1_1$-analysis, $\bf DA$]{schutte1977}. Before presenting the axioms and reduction rules of $\bf SN$, the notion of 
the positive and negative parts of a formula (of $\bf SN$) is in order. 
The notion is not indispensable, but will have the effect for simplifying the subsequent development of $\bf SN$. The notion here to be developed 
following Sch\"utte [1960] \cite{schutteB1960} is somewhat more complicated than its classical counterpart, since we need to distinguish between 
\it antecedent \rm and \it succedent formulas \rm in view of the notion of sequents which is assumed here from the outset. The sequents of 
$\bf SN$ have the form $\Gamma \rightarrow \Delta$ with $\Gamma$ and $\Delta$ being a formula or the empty expression. Unlike the 
traditional Gentzen-type formulation, $\Gamma$ and $\Delta$ represent at most one formula, thus. Since we are implying sequents in place 
of formulas, the positive and negative parts of a formula are defined separately for the antecedent and succedent formulas of a sequent. 
The \it antecedent \rm (\it succedent \rm ) \it positive \rm and \it negative parts \rm of a sequent $A \rightarrow \Delta$ 
($\Gamma \rightarrow B$) is defined recursively as follows: 

\medskip

3.41 $A$ ($B$) is an antecedent (succedent) positive part of $A \rightarrow \Delta$ ($\Gamma \rightarrow B$).

3.42 If $A_1 \wedge A_2$ ($B_1 \vee B_2$) is an antecedent (succedent) positive parts of $\Gamma \rightarrow \Delta$, then $A_1$ and 
$A_2$ ($B_1$ and $B_2$) are antecedent (succedent) positve parts of the sequent.

3.43 $A_1 \vee A_2$ ($B_1 \wedge B_2$) is an antecedent (succedent) negative parts of $\Gamma \rightarrow \Delta$, then $A_1$ and 
$A_2$ ($B_1$ and $B_2$) are antecedent (succedent) negative parts of the sequent.

3.44 $A_1 \supset A_2$ ($B_1 \wedge B_2$) is an antecedent negative parts of $\Gamma \rightarrow \Delta$, then $A_1$ and 
$A_2$ ($B_1$ and $B_2$) are respectively an antecedent positive and an antecedent negative part of the sequent.  

3.45 $\sim A_1$ ($\sim B_1$) is an antecedent (succedent) positive part of $\Gamma \rightarrow \Delta$, then $A_1$ ($B_1$) are antecedent 
(succedent) negative part of the sequent. 

3.46 $\sim A_1$ ($\sim B_1$) is an antecedent (succedent) negative part of $\Gamma \rightarrow \Delta$, then $A_1$ ($B_1$) are antecedent 
(succedent) positive part of the sequent. 

\medskip 

Again, following Sch\"utte, $F[A^+]$ ($F[A^-]$) means that $A$ occurs there as an antecedent positive (negative) part of a sequent, while 
 $G[A_+]$ ($G[A_-]$) menas that $A$ signifies that $A$ is an succedent positive (negative) part of a sequent. Such expressions as 
$F[A^+, B^-] \rightarrow G[C_-, D_-]$ are understood similarly if the specified occrrences of the formulas do not overlap with one another.

For the purpose of illustration, we present some example of antecedent (succedent) positive or negative parts of a sequent. 

\medskip

$F[A^+] \rightarrow G[B_+] = A \rightarrow B,$

$F[A^+] \rightarrow G[B_-] = \enspace \sim \sim A \rightarrow \enspace \sim \sim \sim \sim \sim B,$

$F[A^-, B^+] \rightarrow G[B_-, C_-] = \enspace \sim A \wedge \sim \sim B \rightarrow \enspace \sim C \vee \sim ( B \wedge  \sim A),$

$F[A^+, A^-, B^-] \rightarrow G[A_-, A_+] = \enspace \sim (B \vee (A \supset A)) \rightarrow \enspace A \vee \sim A,$

\medskip 

\noindent where $A, B$ and $C$ are formulas different from each other.

We, further, need the notion of deleting a formulas form another again after Sch\"utte [1960] \cite{schutteB1960}. The notion is dispensed with as 
shown by Sch\"utte [1977] \cite{schutte1977}.

The \it deletion \rm of $A$ from an antecedent formula $F[A^{\pm}]$ (succedent formula $G[A_{\pm}]$), of which the result is expression as 
$F[\; \;^{\pm}]$ ($G[\; \;_{\pm}]$) is defined recursively as follows:

\medskip 

3.51 \enspace If $F[A^+] = A$ ($G[A_+] = A$), then $F[\; \;^+]$ ($G[\; \;_+]$) is the empty expression. 

3.52 \enspace If $F[A^+] = F_1[A \wedge B^+]$ or $= F_1[B \wedge A^+]$ ($G[A_+] = G_1[A \vee B_+]$ or $= G_1[B \vee A_+]$), then $F[\; \; ^+] = F_1[\; \; ^+] \wedge B$ 
($G[\; \; _+] = G_1[\; \; _+] \vee B$).

3.53 \enspace If $F[A^-] = F_1[A \vee B^-]$ or $= F_1[B \vee A^-]$ ($G[A_-] = G_1[A \wedge B_-]$ or $= G_1[B \wedge A_-]$), then 
$F[\; \; ^-] = F_1[\; \; ^-] \wedge \sim B$ ($G[\; \; _-] = G_1[\; \; _-] \vee \sim B$).

3.54 \enspace If $F[A^+] = F_1[A \supset B^-]$, then $F[\; \;^+] = F_1[\; \;^-] \wedge \sim B$.

3.55 \enspace If $F[A^-] = F_1[B \supset A^-]$, then $F[\; \;^+] = F_1[\; \;^-] \wedge B$.

3.56 \enspace If $F[A^+] = F_1[\sim A^-]$, then $F[\; \;^+] = F_1[\; \;^-]$.

3.57 \enspace If $F[A^-] = F_1[\sim A^+]$, then $F[\; \;^-] = F_1[\; \;^+]$.

3.58 \enspace If $G[A_+] = G_1[\sim A_-]$, then $G[\; \;_+] = G_1[\; \;_-]$.

3.59 \enspace If $G[A_-] = G_1[\sim A_+]$, then $G[\; \;_-] = G_1[\; \;_+]$.

\medskip 

An example will be given of the above definition. 

If $F[A^+] \rightarrow G[B_+] = F_1[A \wedge C^+] \rightarrow G_1[B \vee D_+] = A \wedge C \rightarrow B \vee D,$ then 

\noindent $F[\; \;^+] \rightarrow G[\; \;_+] = F_1[\; \;^+] \wedge C \rightarrow G_1[\; \; _+] \vee D = C \rightarrow D.$ 

We, next, present the axioms of $\bf SN$. The formulation below is an adaptation from Sch\"utte [1977] \cite{schutte1977}.

Axiom 1: \enspace $F[A^-] \rightarrow \Delta$ if $A$ is a true constant prime formula.

Axiom 2: \enspace $F[A^+] \rightarrow \Delta$ if $A$ is a false constant prime formula.

Axiom 3: \enspace $\Gamma \rightarrow G[A_+]$ if $A$ is a true constant prime formula.

Axiom 4: \enspace $\Gamma \rightarrow G[A_-]$ if $A$ is a false constant prime formula. 

Axiom 5: \enspace $F[A^+, B^-] \rightarrow \Delta$ if $A$ and $B$ are equivalent formula of length $0$. 

Axiom 6: \enspace $F[A^+] \rightarrow G[B_+] $ if $A$ and $B$ are equivalent formula of length $0$. 

Axiom 7: \enspace $F[A^-] \rightarrow G[B_-] $ if $A$ and $B$ are equivalent formula of length $0$. 

Axiom 8: \enspace $A(x_1, \dots , x_n) \rightarrow B(y_1, \dots , y_m)$ where $n \geq 1$ and $m \geq 0$, or  $n \geq 0$ and 
$m \geq 1$, if for every $ (n + m)$ numerals $k_1, \dots , k_n$, $l_1, \dots , l_m$, $A(k_1, \dots , k_n) \rightarrow B(l_1, \dots , l_m)$ is 
one of the axiom 1-7.


Hereby, the \it length \rm of a formula is defined to be the number of logical symbols. The well-formued formulas of $\bf SN$, 
on the other hand, are defined exactly as in the classical case. The (minimal) parts as indicated in the above axiom 1-7 by $A$ or $B$ 
are called the \it principal parts \rm of these axioms again following Sch\"utte [1977] \cite{schutte1977}. The following reduction rules 
for $\bf SN$ are proposed by Ishimoto [1970 and 197?] \cite{ishimotos, ishimotoe} and Matsuda and A. Ishimoto [1984] \cite{mi} 
(refer to Inou\'e [1984] \cite{inouen} and Shimizu [1990] \cite{shimizu}). The propositional logic part is equivalent to Markov [1950] 
\cite{markov} and Vorob'ev [1952] \cite{voro1}. Probably the system is equivalent to the one proposed by 
Almukdad and Nelson [1984] \cite{an}. 

\bigskip 

The reduction rules of $\bf SN$ are the following. (Below $\Gamma$ and $\Delta$ are only one fomula or empty, respectively.)

\medskip

3.6a 

\medskip 

$\enspace \enspace \displaystyle{\frac{\enspace \enspace F[A \wedge B^-] \rightarrow \Delta }
{\enspace \enspace F[A^-] \rightarrow \Delta \enspace \enspace | \enspace \enspace  F[B^-] \rightarrow \Delta \enspace \enspace  ,}}$ 
$\enspace \enspace (\wedge ^- \rightarrow)$

\bigskip

3.6b 

\medskip  

$\enspace \enspace \displaystyle{\frac{\enspace \enspace \Gamma \rightarrow G[A \wedge B_+] }
{\enspace \enspace \Gamma \rightarrow G[A_+] \enspace \enspace | \enspace \enspace \Gamma \rightarrow G[B_+] \enspace \enspace ,}}$ 
$\enspace \enspace (\rightarrow \wedge_+)$

\bigskip 

3.6c 

\medskip 

$\enspace \enspace \displaystyle{\frac{\enspace \enspace F[A \vee B^+] \rightarrow \Delta }
{\enspace \enspace F[A^+] \rightarrow \Delta \enspace \enspace | \enspace \enspace  F[B^+] \rightarrow \Delta \enspace \enspace  ,}}$ 
$\enspace \enspace (\vee ^+ \rightarrow)$

\bigskip

3.6d 

\medskip  

$\enspace \enspace \displaystyle{\frac{\enspace \enspace \Gamma \rightarrow G[A \vee B_-]}
{\enspace \enspace \Gamma \rightarrow G[A_-] \enspace \enspace | \enspace \enspace \Gamma \rightarrow G[B_-] \enspace \enspace ,}}$ 
$\enspace \enspace (\rightarrow \vee_-)$

\bigskip

3.6e 

\medskip 

$\enspace \enspace \displaystyle{\frac{\enspace \enspace F[A \supset B^+] \rightarrow \Delta }
{\enspace \enspace F[A \supset B^+] \rightarrow \Delta \vee A \enspace \enspace | 
\enspace \enspace  F[A \supset B^+] \wedge B \rightarrow \Delta \enspace \enspace  ,}}$ 
$\enspace \enspace (\supset^+ \rightarrow)$

\bigskip 

3.6f 

\medskip  

$\enspace \enspace \displaystyle{\frac{\enspace \enspace \Gamma \rightarrow G[A \supset B_+] \enspace \enspace  }
{\enspace \enspace \Gamma \wedge A \rightarrow B \enspace \enspace ,}}$ 
$\enspace \enspace (\rightarrow \supset_+)$

\bigskip 

3.6g 

\medskip 

$\enspace \enspace \displaystyle{\frac{\enspace \enspace\Gamma \rightarrow G[A \supset B_-] }
{\enspace \enspace  \Gamma \rightarrow G[\; \;-] \vee A \enspace \enspace | 
\enspace \enspace  \Gamma \rightarrow G[B_-] \enspace \enspace  ,}}$ 
$\enspace \enspace (\supset_- \rightarrow)$

\bigskip 

3.6h

\medskip  

$\enspace \enspace \displaystyle{\frac{\enspace \enspace  \enspace F[\forall xA(x)^+] \rightarrow \Delta \enspace \enspace }
{\enspace \enspace F[\forall xA(x)^+] \wedge A(t) \rightarrow \Delta  \enspace \enspace  ,}}$ 
$\enspace \enspace (\forall ^+ \rightarrow)$

\bigskip 

3.6i

\medskip 

$\enspace \enspace \displaystyle{\frac{\enspace \enspace F[\forall xA(x)^-] \rightarrow \Delta \enspace \enspace}
{\enspace \enspace F[A(b)^-] \rightarrow \Delta   \enspace \enspace  ,}}$
$\enspace \enspace (\forall^- \rightarrow)$

\bigskip 

3.6j 

\medskip  

$\enspace \enspace \displaystyle{\frac{\enspace \enspace  \enspace \Gamma \rightarrow G[\forall xA(x)_+] \enspace \enspace }
{\enspace \enspace \Gamma \rightarrow A(b)  \enspace \enspace  ,}}$ 
$\enspace \enspace (\rightarrow \forall_+)$

\bigskip 

3.6k 

\medskip  

$\enspace \enspace \displaystyle{\frac{\enspace \enspace  \enspace \Gamma \rightarrow G[\forall xA(x)_-] \enspace \enspace }
{\enspace \enspace \Gamma \rightarrow G[\forall xA(x)_-]  \vee \sim A(t)  \enspace \enspace  ,}}$ 
$\enspace \enspace (\rightarrow \forall_-)$

\bigskip

3.6l 

\medskip  

$\enspace \enspace \displaystyle{\frac{\enspace \enspace  \enspace F[\exists xA(x)^+] \rightarrow \Delta \enspace \enspace }
{\enspace \enspace F[A(b)^+] \rightarrow \Delta  \enspace \enspace  ,}}$ 
$\enspace \enspace (\exists ^+ \rightarrow)$

\bigskip 

3.6m 

\medskip 

$\enspace \enspace \displaystyle{\frac{\enspace \enspace F[\exists xA(x)^-] \rightarrow \Delta \enspace \enspace}
{\enspace \enspace F[\exists xA(x)^-] \wedge \sim A(t) \rightarrow \Delta   \enspace \enspace  ,}}$
$\enspace \enspace (\exists^- \rightarrow)$

\bigskip 

3.6n 

\medskip  

$\enspace \enspace \displaystyle{\frac{\enspace \enspace  \enspace \Gamma \rightarrow G[\exists xA(x)_+] \enspace \enspace }
{\enspace \enspace \Gamma \rightarrow  G[\exists xA(x)_+] \vee A(t) \enspace \enspace  ,}}$ 
$\enspace \enspace (\rightarrow \exists_+)$

\bigskip 

3.6o 

\medskip  

$\enspace \enspace \displaystyle{\frac{\enspace \enspace  \enspace \Gamma \rightarrow G[\exists xA(x)_-] \enspace \enspace }
{\enspace \enspace \Gamma \rightarrow \sim A(b) \enspace \enspace  ,}}$ 
$\enspace \enspace (\rightarrow \exists_-)$


\bigskip 

\noindent where the $b$'s in $(\exists ^+ \rightarrow)$, $(\rightarrow \forall_+)$, $(\exists ^+ \rightarrow)$ and $(\rightarrow \exists_-)$ 
are variables not occurring free in sequent to be reduced, i.e., proper (eigen) variables, and $t$ is a term. (In the sequent, we will often 
use the words a \it principal formula \rm and a \it quasi-principal formula \rm of a reduction as usual.) A sequent 
$\Gamma \rightarrow \Delta$ of $\bf SN$ is provable in $\bf SN$ if each branch of the tableau obtained from $\Gamma \rightarrow \Delta$ by 
applying reduction rules 3.6a--3.6o ends with one of the axioms of $\bf SN$.  $\bf SN \it \vdash \Gamma \rightarrow \Delta$ means that 
it is provable in $\bf SN$. A formula of $A$ of $\bf SN$ is a \it theorem \rm of $\bf SN$ if $\bf SN \vdash \rightarrow A$. We also write 
$\bf SN \vdash A$ if $A$ is a theorem of $\bf SN$ in the case where no ambiguity arises. 

We have to mention the notion of \it normal tableaus \rm for $\bf SN$. In order to prove the completeness of $\bf SN$, we need 
normal tableaus with the notion of Hintikka formula. However, this paper is not intended to develop semantical study. So in this paper,
 we shall not go further on this direction. 

For a semantic study for the simpler predicate logic part 
of $\bf SN$, refer to Shimizu [1990] \cite{shimizu}. The predicate logic part of $\bf SN$ is due to Prof. Arata Ishimoto. It is also used in 
Matsuda and Ishimoto [1984] \cite{mi} and Inou\'e [1984, 202?] \cite{inoueac, inouefunda} for logic programming and Prolog. 

\section{Subsystems $\bf PCN$ and $\bf FN$ of $\bf SN$}

Now, we propose two subsystems $\bf FN$ 
and $\bf PCN$ of $\bf SN$. The first subsystem $\bf FN$ is obtained from $\bf SN$ by supressing the reduction rule $(\supset^+ \rightarrow)$. 
The second one $\bf PCN$ may be called a \it pseudo-classical subsystem \rm of $\bf SN$. The axioms and reduction rules of $\bf PCN$ are 
presented as follows:

\medskip

(Axioms of $\bf PCN$)

\medskip 

3.71 \enspace \enspace $F[A^-] \rightarrow \enspace $ if $A$ is a true constant prime formula. 

3.72 \enspace \enspace $F[A^+] \rightarrow \enspace $ if $A$ is a false constant prime formula. 

3.73 \enspace \enspace $F[A^+, B^-] \rightarrow \enspace $ if $A$ and $B$ are equivalent formulas of length $0$.

3.74 \enspace \enspace $A(x_1 , \dots , x_n) \rightarrow \enspace $ where $n \geq 1$ and for every numerals $k_1, \dots , k_n$ 
$A(k_1, \dots , k_n) \rightarrow \enspace $ is one of the axioms 3.71, 3.72 ad 3.73. 

\medskip 

(Reduction rules of $\bf PCN$)

\medskip 

3.81 

\medskip 

$\enspace \enspace \displaystyle{\frac{\enspace \enspace F[A \wedge B^-] \rightarrow \enspace \enspace }
{\enspace \enspace F[A^-] \rightarrow  \enspace \enspace | \enspace \enspace  F[B^-] \rightarrow \enspace \enspace  ,}}$ 
$\enspace \enspace (\wedge ^- \underset{c}{\rightarrow})$

\bigskip 

3.82  

\medskip 

$\enspace \enspace \displaystyle{\frac{\enspace \enspace F[A \vee B^+] \rightarrow \enspace \enspace }
{\enspace \enspace F[A^+] \rightarrow \enspace \enspace | \enspace \enspace  F[B^+] \rightarrow  \enspace \enspace  ,}}$ 
$\enspace \enspace (\vee ^+  \underset{c}{\rightarrow})$

\bigskip

3.83 

\medskip  

$\enspace \enspace \displaystyle{\frac{\enspace \enspace  \enspace F[\forall xA(x)^+] \rightarrow  \enspace \enspace }
{\enspace \enspace F[\forall xA(x)^+] \wedge A(t) \rightarrow  \enspace \enspace  ,}}$ 
$\enspace \enspace (\forall ^+  \underset{c}{\rightarrow})$

\bigskip 

3.84 

\medskip 

$\enspace \enspace \displaystyle{\frac{\enspace \enspace F[\forall xA(x)^-] \rightarrow \enspace \enspace}
{\enspace \enspace F[A(b)^-] \rightarrow \enspace \enspace  ,}}$
$\enspace \enspace (\forall^- \underset{c}{\rightarrow})$

\bigskip 

3.85  

\medskip  

$\enspace \enspace \displaystyle{\frac{\enspace \enspace  \enspace F[\exists xA(x)^+] \rightarrow \enspace \enspace }
{\enspace \enspace F[A(b)^+] \rightarrow  \enspace \enspace  ,}}$ 
$\enspace \enspace (\exists ^+  \underset{c}{\rightarrow})$

\bigskip 

3.86  

\medskip 

$\enspace \enspace \displaystyle{\frac{\enspace \enspace F[\exists xA(x)^-] \rightarrow  \enspace \enspace}
{\enspace \enspace F[\exists xA(x)^-] \wedge \sim A(t) \rightarrow \enspace \enspace ,}}$
$\enspace \enspace (\exists^-  \underset{c}{\rightarrow})$

\bigskip 

\noindent where $b$'s are proper (eigen) variables and $t$ is a term. 

We shall continue to use the name of axioms and reduction rules of $\bf SN$ in place of the above ones from now on because axioms 3.71--3.74 
and 3.81--3.86 are their special cases of their conterprts in $\bf SN$ with the lack of sucfcedent formula. 

For example, we will use the Axiom 8 in place of the Axiom 3.74, the reduction rule  $(\wedge ^- \rightarrow)$ 
inplace of $(\wedge ^- \underset{c}{\rightarrow})$. It is remarked that $\bf PCN$ is classical in its appearance as long as we confine ourselves 
to antecedent formulas. A formula of $\bf PCN$ may contain implication sign though reduction rules for implication are not available in $\bf PCN$. 
The point will be crucial in the sequel to show that the cut elimination theorem is provable with in $\bf PCN$. The fact will lead us to 
a proof of the embedding theorem. 

The subsystems $\bf FN$ and $\bf PCN$ here introduced are very special systems.  It is, however, believed that these 
logics are worthy of more careful studies. 

\begin{proposition} \mbox{} \\ 

$(1)$ \enspace \enspace $\bf FN \it \vdash \Gamma \rightarrow \Delta \enspace \Rightarrow \enspace 
\bf SN \it \vdash \Gamma \rightarrow \Delta$.

$(2)$ \enspace \enspace $\bf PCN \it \vdash \Gamma \rightarrow \enspace \Rightarrow \enspace \bf SN \it \vdash \Gamma \rightarrow$.
\end{proposition}

\it Proof. \rm Trivial from the definitions. $\Box$


\section{Theorems and metatheorems of $\bf SN$}

We wish to obtain a number of theorems and metatheorems of $\bf SN$ which will be in order for proving main theorems to be stated \S 5 and \S 6. 

\begin{theorem} \rm (Simultaneous substitution theorem) \it 

$$\bf SN \it \vdash \Gamma \rightarrow \Delta \enspace \Rightarrow \enspace 
\bf SN \it \vdash (\Gamma \rightarrow \Delta)[t_1/x_1, \dots ,t_n/x_n],$$

\noindent where $\Gamma \rightarrow \Delta)[t_1/x_1, \dots ,t_n/x_n]$ stands for the result of simultaneously replacing $x_1, \dots , X_n$ in 
$\Gamma \rightarrow \Delta$ by terms $t_1, \dots , t_n$, respectively, if variables occurring free in $\Gamma \rightarrow \Delta$ are among 
$x_1 , \dots , X_n$ and no terms $t_1, \dots , t_n$ are variables occurring bound in $\Gamma \rightarrow \Delta$ and contain any variables 
occurring bound in $\Gamma \rightarrow \Delta$.
\end{theorem}

\it Proof. \rm This is easily proved by induction on the length of the tableau.  $\Box$ 

\begin{theorem} \rm (Inversion theorems) \it 

\medskip 

$(1)$ \enspace \enspace $\bf SN \it \vdash F[A \wedge B^-] \rightarrow \Delta \enspace \Leftrightarrow \enspace 
(\bf SN \it \vdash F[A^-] \rightarrow \Delta \mbox{ and } \bf SN \it \vdash F[B^-] \rightarrow \Delta)$.

$(2)$ \enspace \enspace $\bf SN \it \vdash \Gamma \rightarrow G[A \wedge B_+] \enspace \Leftrightarrow \enspace 
(\bf SN \it \vdash \Gamma \rightarrow G[A_+] \mbox{ and  } \bf SN \it \vdash \Gamma \rightarrow G[B_+])$.

$(3)$ \enspace \enspace $\bf SN \it \vdash F[A \vee B^+] \rightarrow \Delta \enspace \Leftrightarrow \enspace 
(\bf SN \it \vdash F[A^+] \rightarrow \Delta \mbox{ and } \bf SN \it \vdash F[B^+] \rightarrow \Delta)$.

$(4)$ \enspace \enspace $\bf SN \it \vdash \Gamma \rightarrow G[A \vee B_-] \enspace \Leftrightarrow \enspace 
(\bf SN \it \vdash \Gamma \rightarrow G[A_-] \mbox{ and  } \bf SN \it \vdash \Gamma \rightarrow G[B_-])$.

$(5)$ \enspace \enspace $\bf SN \it \vdash \Gamma \rightarrow G[A \supset B_-] \enspace \Leftrightarrow \enspace 
(\bf SN \it \vdash \Gamma \rightarrow G[\; \;_-] \vee A \mbox{ and  } \bf SN \it \vdash \Gamma \rightarrow G[B_-])$.

$(6)$ \enspace \enspace $\bf SN \it \vdash F[\forall xA(x)^-] \rightarrow \Delta \enspace \Leftrightarrow \enspace 
\bf SN \it \vdash F[A(b)^-] \rightarrow \Delta $.

$(7)$ \enspace \enspace $\bf SN \it \vdash \Gamma \rightarrow G[\forall xA(x)_+] \enspace \Leftrightarrow \enspace 
\bf SN \it \vdash \Gamma \rightarrow A(b)$.

$(8)$ \enspace \enspace $\bf SN \it \vdash F[\exists xA(x)^+] \rightarrow \Delta \enspace \Leftrightarrow \enspace 
\bf SN \it \vdash F[A(b)^+] \rightarrow \Delta $.

$(9)$ \enspace \enspace $\bf SN \it \vdash \Gamma \rightarrow G[\exists xA(x)_-] \enspace \Leftrightarrow \enspace 
\bf SN \it \vdash \Gamma \rightarrow \sim A(b)$.

\end{theorem}

\it Proof. \rm We only prove $\Rightarrow$ of (1) on the basis of the given sequent, contain the principal formula $A \wedge B$ of 
the reduction rule in the tableau for the given sequent $F[A \wedge B^-] \rightarrow \Delta$. 

\it Basis \rm Case 1: let $F[A \wedge B^-] \rightarrow \Delta$ be an axiom.  $F[A^-] \rightarrow \Delta$ and 
$F[B^-] \rightarrow \Delta$ are also axioms, since the antecedent negative part $A \wedge B$ of 
$F[A \wedge B^-] \rightarrow \Delta$ is by no means a principal part of the axiom.

Case 2: Suppose that  $F[A \wedge B^-] \rightarrow \Delta$ is reduced by applying $(\wedge ^- \rightarrow)$ in the following way: 

\medskip 

$\enspace \enspace \displaystyle{\frac{\enspace \enspace F[A \wedge B^-] \rightarrow \Delta }
{\enspace \enspace F[A^-] \rightarrow \Delta \enspace \enspace | \enspace \enspace  F[B^-] \rightarrow \Delta \enspace \enspace  .}}$ 
$\enspace \enspace (\wedge ^- \rightarrow)$

\bigskip

\noindent It, then, follows that $F[A^-] \rightarrow \Delta$ and  $F[B^-] \rightarrow \Delta$ are both provable in $\bf SN$. 

\it Induction steps \rm Case 1: Assume that $F[A \wedge B^-] \rightarrow \Delta$ is reduced by appplying reduction rules of which 
the principal formula is not $A \wedge B$. For example, let us asssume that the reduction rule applied to the given sequent be 
$(\vee^+ \rightarrow)$: 

\medskip 

$\enspace \enspace \displaystyle{\frac{\enspace \enspace F_1[A \wedge B^-, C \vee D^+] \rightarrow \Delta }
{\enspace \enspace F_1[A \wedge B^-, C^+] \rightarrow \Delta \enspace \enspace | \enspace \enspace  
F_1[A \wedge B^-, D^+] \rightarrow \Delta \enspace \enspace  ,}}$ 
$\enspace \enspace (\vee^+ \rightarrow)$

\bigskip

\noindent where $F_1[A \wedge B^-, C \vee D^+] \rightarrow \Delta =  F[A \wedge B^-] \rightarrow \Delta$.  These two sequents 
obtained as  a result of the reduction are, then, subject to the hypothesis of induction, since the corresponding to the given sequent. 
These sequents are provable in $\bf SN$ by the hypothesis of induction. We, thus, obtain 
$\bf SN \it \vdash F_1[A^-, C \vee D^+] \rightarrow \Delta$ 
and $\bf SN \it \vdash F_1[B^-, C \vee D^+] \rightarrow \Delta$. 

Case 2: Assume that the given sequent is reduced by applying $(\wedge^- \rightarrow)$ and that the principal formula of the rule is 
$A \wedge B$. The case is similarly proved as in the case 2 of the basis case.

In the cases that the rules to be inverted be succedent rules, there are more basis cases than in the antededent cases, but this 
does not present any difficulties. In addition, it is remarked that the tableau obtained as a result of inversion is not longer 
than the original one. 

The converse, namely, $\Leftarrow$ holds obviously, since $F[A \wedge B^-] \rightarrow \Delta$ is reduced to $F[A^-] \rightarrow \Delta$ 
and $F[B^-] \rightarrow \Delta$ by applying $(\wedge^- \rightarrow)$.  $\Box$

\begin{theorem} \rm (Contraction theorems) \it 

\medskip 

$(1)$ \enspace \enspace $\bf SN \it \vdash F[A^\pm, A^\pm] \rightarrow \Delta \enspace \Rightarrow \enspace 
\bf SN \it \vdash F[A^\pm , \; \;^\pm ]  \rightarrow \Delta$.

$(2)$ \enspace \enspace $\bf SN \it \vdash \Gamma \rightarrow G[A_\pm, A_\pm] \enspace \Rightarrow 
\enspace \bf SN \it \vdash \Gamma \rightarrow G[A_\pm, \; \;_\pm]$.

\end{theorem}

This is proved by induction on the length of the tableau. In the proof, the inversion theorems are indispensable for the treatment of the 
induction steps. Use is made of the cace that the tableau of any sequent obtained by applying the inversion theorem is not longer 
than the original one. In addition, it is noticed that the presence of the quasi-principal formula for five rules 
$(\supset^+ \rightarrow)$, $(\forall^+ \rightarrow)$, $(\rightarrow \forall_+)$, $(\exists^-\rightarrow)$ and $(\rightarrow \exists_+)$ is 
required in view of the failure of the inversion theorems for them.

We are, now, proceeding to metatheorems of $\bf SN$ corresponding to structural rules. In what follows, those theorems of which we do not
 give proofs are demonstrated with ease by induction on the length of the tableau. 

\begin{theorem} \rm (Thinning theorems) \it 

\medskip 

$(1)$ \enspace \enspace $\bf SN \it \vdash F[\; \;^\pm] \rightarrow \Delta \enspace \Rightarrow \enspace 
\bf SN \it \vdash F[A^\pm] \rightarrow \Delta$.

$(2)$ \enspace \enspace $\bf SN \it \vdash \Gamma \rightarrow G[\; \;_\pm] \enspace \Rightarrow 
\enspace \bf SN \it \vdash \Gamma \rightarrow G[A_\pm]$.
\end{theorem} 

\begin{theorem} \rm (Translation theorems) \it 

\medskip 

$(1)$ \enspace \enspace $\bf SN \it \vdash F[A^\pm, \; \;^\pm] \rightarrow \Delta \enspace \Rightarrow \enspace 
\bf SN \it \vdash F[\; \;^\pm , A^\pm ]  \rightarrow \Delta$.

$(2)$ \enspace \enspace $\bf SN \it \vdash \Gamma \rightarrow G[A_\pm, \; \;_\pm] \enspace \Rightarrow 
\enspace \bf SN \it \vdash \Gamma \rightarrow G[\; \;_\pm, A_\pm]$.
\end{theorem}

\it Proof\rm .  We prove (1) only. The proof of (2) is similar to that of (1). Assume  $\bf SN \it \vdash F[A^\pm, \; \;^\pm] \rightarrow \Delta$. 
By the thinning theorem,  $\bf SN \it \vdash F[A^\pm, A^\pm] \rightarrow \Delta$. We obtain
 $\bf SN \it \vdash F[\; \;^\pm, A^\pm] \rightarrow \Delta$ by the constraction theorem.    $\Box$

\begin{theorem} \rm (Interchange theorems) \it 

\medskip 

$(1)$ \enspace \enspace $\bf SN \it \vdash F[A^\pm, B^\pm] \rightarrow \Delta \enspace \Rightarrow \enspace 
\bf SN \it \vdash F[B^\pm , A^\pm ]  \rightarrow \Delta$.

$(2)$ \enspace \enspace $\bf SN \it \vdash \Gamma \rightarrow G[A_\pm, B_\pm] \enspace \Rightarrow 
\enspace \bf SN \it \vdash \Gamma \rightarrow G[B_\pm, A_\pm]$.
\end{theorem}

\it Proof\rm . 

\smallskip 

\enspace \enspace $\bf SN \it \vdash F[A^+, B^+] \rightarrow \Delta$

$\Rightarrow \bf SN \it \vdash F[\; \;^+, B^+] \wedge A \rightarrow \Delta$ \enspace \enspace  (Translation theorem)

$\Rightarrow \bf SN \it \vdash F[B^+, A^+] \rightarrow \Delta .$ \enspace \enspace  (Translation theorem)

\medskip 

The other cases are analogously dealt with. $\Box$

We shall present the disjunction and existence properties of $\bf SN$ as rather more general forms.

\begin{theorem} \rm (Disjunction property of $\bf SN$) \it 
$$\mbox{\bf SN} \vdash \enspace \rightarrow G[A_+] \enspace \Leftrightarrow \enspace 
(\mbox{\bf SN}  \vdash \enspace \rightarrow G[\; \;_+] \mbox{ or } \enspace \mbox{\bf SN} \vdash \enspace \rightarrow A).$$
\end{theorem}

\it Proof of $\Rightarrow$ for Theorem 5.7\rm .  If $\rightarrow G[\; \;_+]$ is empty, then $\mbox{\bf SN} \vdash \enspace \rightarrow A$ is 
trivial by 3.51. So assume that $\rightarrow G[\; \;_+]$ is not empty. We prove the theorem by induction on the number of 
the sequents which are successively obtained by reducing the given 
$\rightarrow G[A_+]$ down to that where $A$ occurs for the first time.  

\it Basis \rm Case 1: Let $\rightarrow G[A_+]$ be an axiom. If $A$ does not contain the principal part of the axiom, then $\rightarrow G[\; \;_+]$ 
is also an axiom. Otherwise  $\rightarrow A$ constitutes an aaxiom.

Case 2: Let $G[A_+]$ be $G_1[A_+, B \supset C_+]$ with $B \supset C$ being the principal formula. The sequent, then, is reduced as follows:

\medskip  

$\enspace \enspace \displaystyle{\frac{\enspace \enspace \rightarrow G_1[A_+, B \supset C_+] \enspace \enspace  }
{\enspace \enspace B \rightarrow C \enspace \enspace .}}$ 
$\enspace \enspace (\rightarrow \supset_+)$

\bigskip 

\noindent We, then, obtain $\bf SN \it \vdash \enspace \rightarrow G_1[\; \;_+, B \supset C_+]$, i.e., $\rightarrow G[\; \;_+]$ which is forthcoming by:

\medskip  

$\enspace \enspace \displaystyle{\frac{\enspace \enspace \rightarrow G_1[\; \;_+, B \supset C_+] \enspace \enspace  }
{\enspace \enspace B \rightarrow C \enspace \enspace .}}$ 
$\enspace \enspace (\rightarrow \supset_+)$

\bigskip 

Case 3: Let $G[A_+]$ be $G_1[B \supset C_+]$, where $B \supset C$ constitutes a positive part of $A$. The given tableau, then is of the form: 

\medskip  

$\enspace \enspace \displaystyle{\frac{\enspace \enspace \rightarrow G_1[\; \;_+, B \supset C_+] \enspace \enspace  }
{\enspace \enspace B \rightarrow C \enspace \enspace .}}$ 
$\enspace \enspace (\rightarrow \supset_+)$

\bigskip 

\noindent It is clear that  $\rightarrow A$ is provable in $\bf SN$, since $\rightarrow A$ is reduced to $B \rightarrow C$ by applying 
$(\rightarrow \supset_+)$. 

Case 4: Let $G[A_+]$ be $G_1[A_+, \forall xB(x)_+]$, and consider the following reduction: 

\medskip  

$\enspace \enspace \displaystyle{\frac{\enspace \enspace  \rightarrow G_1[A_+,\forall xB(x)_+] \enspace \enspace  }
{\enspace \enspace \rightarrow B(b) \enspace \enspace .}}$ 
$\enspace \enspace (\rightarrow \forall_+)$

\bigskip 

We obtain $\bf SN \it \vdash  \rightarrow G_1[\; \;_+,\forall xB(x)_+] $, since the sequent is reduced to $\rightarrow B(b)$ by applying the 
same rule. 

Case 5: Let $G[A_+]$ be $G_1[\forall xB(x)_+]$. Assume, further, that $A$ contains $\forall xB(x)$ as its positive part. 
If $\rightarrow G_1[\forall xB(x)_+]$ has the following reduction: 

\medskip  

$\enspace \enspace \displaystyle{\frac{\enspace \enspace  \rightarrow G_1[\forall xB(x)_+] \enspace \enspace  }
{\enspace \enspace \rightarrow B(b) \enspace \enspace .}}$ 
$\enspace \enspace (\rightarrow \forall_+)$

\bigskip 

\noindent Then, $\rightarrow A$ is provable in $\bf SN$, since the sequent is reduced to $\rightarrow B(b)$ by applying $(\rightarrow \forall_+)$. 

Case 6: Let $G[A_+]$ be $G_1[\exists xB(x)_+]$, and consider that $\rightarrow G[A_+]$ is reduced to $\rightarrow \enspace \sim B(b)$ 
by applying  $(\rightarrow \exists_+)$.  Although the case splits up into two subcases like the case 4 and 5 of the baseis case, they are similarly 
dealt with as in the cases 4 and 5.

\it Induction steps \rm Case 1: Let $G[A_+]$ be $G_1[A_+, B \wedge C_+]$, which results by way of $(\rightarrow \wedge_+)$ in the following way: 

\medskip  

$\enspace \enspace \displaystyle{\frac{\enspace \enspace \rightarrow G[A_+, B \wedge C_+]}
{\enspace \enspace \rightarrow G[A_+, B_+] \enspace \enspace | \enspace \enspace \rightarrow G[A_+, C_+] \enspace \enspace .}}$ 
$\enspace \enspace (\rightarrow \wedge_+)$

\bigskip 

By H.I. (we shall abbreviate the hypothesis of induction as H.I.), we obtain ($\bf SN \it \vdash \rightarrow A$ or 
$\bf SN \it \vdash \rightarrow G_1[\; \;_+, B_+]$) and ($\bf SN \it \vdash \rightarrow A$ or $\bf SN \it \vdash \rightarrow G_1[\; \;_+, C_+]$). 
If  $\rightarrow A$ is not provable in $\bf SN$, then $\bf SN \it \vdash \rightarrow G_1[\; \;_+, B_+]$ and $\bf SN \it \vdash \rightarrow G_1[\; \;_+, C_+]$. 
From this we obtain $\bf SN \it \vdash \rightarrow G[\; \;_+, B \wedge C_+]$ by $(\rightarrow \wedge_+)$. 

Case 2:  Let $G[A_+]$ be of the form $G[A[B \wedge C_+]_+]$ with $B \wedge C$ occurring in $A$ as a succedent positive part, consequently, 
as a positive part of $G$, It is, further, assumed that the given sequent be reduced in the following way:

\medskip  

$\enspace \enspace \displaystyle{\frac{\enspace \enspace \rightarrow G[A[B \wedge C_+] _+}
{\enspace \enspace \rightarrow G[A[B_+] _+] \enspace \enspace | \enspace \enspace \rightarrow G[A[C_+]_+] \enspace \enspace .}}$ 
$\enspace \enspace (\rightarrow \wedge_+)$

\bigskip 

By H.I., we have ($\bf SN \it \vdash \rightarrow A[B_+]$ or 
$\bf SN \it \vdash \rightarrow G[\; \;_+]$) and ($\bf SN \it \vdash \rightarrow A[C_+]$ or $\bf SN \it \vdash \rightarrow G[\; \;_+]$). If 
$\bf SN \it \vdash \rightarrow G[\; \;_+]$, we have done it. 
Otherwise, we have $\bf SN \it \vdash \rightarrow A[B_+]$ and $\bf SN \it \vdash \rightarrow A[C_+]$ to which 
$\bf SN \it \vdash \rightarrow A[B \wedge C]_+]$ is reduced by way of $(\rightarrow \wedge_+)$.

Case 3: Assume that $\bf SN \it \vdash \rightarrow G[A_+]$ is reduced by applying $(\rightarrow \vee_-)$ and that the principal formula of
 the reduction is $B \vee C$. The case splits up into two subcases in accordance with $B \vee C$ being or not being contained in $A$ as 
its negative part. These are taken care of  analogously with the cases 1 and 2 of the induction steps. 

Case 4: Let $G[A_+]$ be $G_1[A_+, B \supset C_-]$. Suppose that $\rightarrow G_1[\; \;_+, B \supset C_-]$ is reduced as follows:

\medskip  

$\enspace \enspace \displaystyle{\frac{\enspace \enspace \rightarrow G[A_+, B \supset C_-] }
{\enspace \enspace \rightarrow G[A_+, \; \;_-] \vee B \enspace \enspace | \enspace \enspace \rightarrow G[A_+, C_-] \enspace \enspace .}}$ 
$\enspace \enspace (\rightarrow \supset_-)$

\bigskip 

We, then, have ($\bf SN \it \vdash \rightarrow A$ or 
$\bf SN \it \vdash \rightarrow G[\; \;_+, \; \;_-] \vee B$) and ($\bf SN \it \vdash \rightarrow A$ or $\bf SN \it \vdash \rightarrow G[\; \;_+, C_-]$) 
by H.I.. If $\rightarrow A$ is provable, we are satisfied with it. 

In the contrary case, $\rightarrow G[\; \;_+, \; \;_-] \vee B$ and $\rightarrow G[\; \;_+, C_-]$ are provable in $\bf SN$, and 
$\rightarrow G[A_+, B \supset C_-]$ is 
reduced to them by $(\rightarrow \supset_-)$. 

Case 5:  Let $G[A_+]$ be of the form $G[A[B \supset C_-]_+]$ with $B \supset C$ occurring in $A$ as a succedent negative part, consequently, 
as a negative part of $G$. It is also assumed that the given sequent be reduced in the following way:

\medskip  

$\enspace \enspace \displaystyle{\frac{\enspace \enspace \rightarrow G[A[B \supset C_-] _+]}
{\enspace \enspace \rightarrow G[A[\; \;_-] _+] \vee B \enspace \enspace | \enspace \enspace \rightarrow G[A[C_+]_+] \enspace \enspace .}}$ 
$\enspace \enspace (\rightarrow \supset_-)$

\bigskip 

By the translation theorem, $\rightarrow G[A[\; \;_-]_+] \vee B$ gives rise to $\rightarrow G[A[\; \;_-]_+, \vee B]$, which is proved by a tableau not 
longer than that for the former. By H.I., we have ($\bf SN \it \vdash \rightarrow A[\; \;_-] \vee B$ or 
$\bf SN \it \vdash \rightarrow G[\; \;_+]$) and ($\bf SN \it \vdash \rightarrow A[C_-]$ or $\bf SN \it \vdash \rightarrow G[\; \;_+]$). 
If $\bf SN \it \vdash \rightarrow G[\; \;_+]$, we have done it. Otherwise, $\rightarrow G[A[B \supset C_-]_+]$ is reduced to 
$\rightarrow A[\; \;_-] \vee B$ and $\rightarrow A[C_-]$  by way of $(\rightarrow \supset_-)$. 

Case 6: Let $G[A_+]$ be $G_1[A_+, \forall xB(x)_-]$. Suppose that the given sequent is reduced by the following:

\medskip  

$\enspace \enspace \displaystyle{\frac{\enspace \enspace  \rightarrow G_1[A_+,\forall xB(x)_-] \enspace \enspace  }
{\enspace \enspace \rightarrow G_1[A_+, \forall xB(x)_-] \vee \sim B(t) \enspace \enspace .}}$ 
$\enspace \enspace (\rightarrow \forall_-)$

\bigskip 

By H.I., we obtain $\bf SN \it \vdash  \rightarrow A$ or $\bf SN \it \vdash  \rightarrow G_1[\; \;_+,\forall xB(x)_-] \vee \sim B(t)$. 
If $\bf SN \it \vdash  \rightarrow A$ is not provable, then 
$\bf SN \it \vdash  \rightarrow G_1[\; \;_+,\forall xB(x)_-] \vee \sim B(t)$ holds. It follows that 
$\bf SN \it \vdash  \rightarrow G_1[\; \;_+,\forall xB(x)_-]$. The tableau thereof is of the form: 

\medskip  

$\enspace \enspace \displaystyle{\frac{\enspace \enspace  \rightarrow G_1[\; \;_+,\forall xB(x)_-] \enspace \enspace  }
{\enspace \enspace \rightarrow G_1[\; \;_+, \forall xB(x)_-] \vee \sim B(t) \enspace \enspace .}}$ 
$\enspace \enspace (\rightarrow \forall_-)$

\bigskip 

Case 7: Let $G[A_+]$ be $G_1[\forall xB(x)_-]$. Assume, further, that $A$ contains $\forall xB(x)$ as its negative part and that  $G[A_+]$ is 
reduced by applying $(\rightarrow \forall_-)$, and the principal formula of the reduction is $\forall xB(x)$. We can take care of the case as in 
the case 6. Note that we use the translation theorem as in the case 5 of the induction steps.

Case 8: Let $G[A_+]$ be $G_1[\exists xB(x)_+]$. Suppose that $A$ contains $\exists xB(x)$ as its positive part and that $\rightarrow G[A_+]$ 
is reduced by $(\rightarrow \exists_+)$ with the principal formula of the reduction being $\exists xB(x)$. The case is treated analogously 
in the cases 6 and 7. 

\it Proof of $\Leftarrow$ for Theorem 5.7\rm . If is immediate that  $\bf SN \it \vdash  \rightarrow G_1[A_+]$ holds by means of the 
thinning theorem and the translation theorem. $\Box$

\begin{theorem} \rm (E-theorem of $\bf SN$) \it 
$$\mbox{\bf SN} \vdash \enspace \rightarrow G[\exists xA(x)_+] \enspace \Leftrightarrow \enspace 
(\mbox{A term $t$ is found such that \bf SN} \vdash \enspace \rightarrow G[A(t)_+]).$$
\end{theorem}

\it Proof of $\Rightarrow$ for Theorem 5.8\rm . We prove the theorem by induction on the length of the given tableau for $\rightarrow G[\exists xA(x)_+]$
 down to the sequents where $\rightarrow G[\exists xA(x)_+]$ was first introduced. (In what follow the length will be understood in this 
generalezed sense.)

\it Basis \rm Case 1: Assume the given sequent is an axiom. In view of the definition of axioms, 
$\exists xA(x)$ is not the principal part of an axiom. Thus, $\rightarrow G[A(t)_+]$ is 
also an axiom for any term $t$ of $\bf SN$.

Case 2: Assume that $G[\exists xA(x)_+]$ of the given sequent is introduced by thinning as a result of reduction rules $(\rightarrow \supset_+)$, 
$(\rightarrow \forall_+)$ or $(\rightarrow \exists_-)$. For example, let $G[\exists xA(x)_+]$ be $G_1[\exists  xA(x)_+, B \supset C_+]$. Suppose 
the given sequent is subject to the following reduction:

\medskip  

$\enspace \enspace \displaystyle{\frac{\enspace \enspace \rightarrow G_1[\exists xA(x)_+, B \supset C_+] \enspace \enspace  }
{\enspace \enspace B \rightarrow C \enspace \enspace .}}$ 
$\enspace \enspace (\rightarrow \supset_+)$

\bigskip 

\noindent It, then, follows that, for any term $t$, $\rightarrow G_1[A(t)_+, B \supset C_+]$. The sequent is reduced to $B \rightarrow C$ 
by applying the same rule again. Other cases are similarly dealt with.

\it Induction steps \rm We shall confine ourselves to three typical cases. The cases for $(\rightarrow \vee_-)$, 
$(\rightarrow \supset_-)$ or $(\rightarrow \forall_-)$ are proved in the analogous way.

Let $G[\exists xA(x)_+]$ be $G_1[\exists xA(x)_+, B \wedge C_+]$. Assume the given sequent is reduced in the form:

\medskip  

$\enspace \enspace \displaystyle{\frac{\enspace \enspace \rightarrow G_1[\exists xA(x)_+, B \wedge C_+]}
{\enspace \enspace \rightarrow G[\exists xA(x)_+, B_+] \enspace \enspace | \enspace  \rightarrow G[\exists xA(x)_+, C_+] \enspace \enspace .}}$ 
$\enspace \enspace (\rightarrow \wedge_+)$

\bigskip 

\noindent By Theorem 5.7, we obtain ($\bf SN \it \vdash  \rightarrow \exists xA(x)$ or $\bf SN \it \vdash  \rightarrow G_1[\; \;_+, B_+]$) and 
($\bf SN \it \vdash  \rightarrow \exists xA(x)$ or $\bf SN \it \vdash  \rightarrow G_1[\; \;_+, C_+]$). Assume that $\rightarrow \exists xA(x)$ 
is provable. A term $t$ is, then, found by H..I. such that $\bf SN \it \vdash A(t)$, since tableau for $\exists xA(x)$ is not longer than that for 
$\rightarrow G_1[\exists xA(x)_+, B \wedge C_+]$ as noticed in the proof of Theorem 5.7. From this follows $\rightarrow G[A(t)_+]$ 
by thinning. In case where $\rightarrow G_1[\; \;_+, B_+]$ and $\rightarrow G_1[\; \;_+, C_+]$, we obtain $\rightarrow G[A(t)_+, B\wedge C_+]$ 
by $(\rightarrow \wedge_+)$ and thinning.

Let $G[\exists xA(x)_+]$ be $G_1[\exists xA(x)_+, \exists xB(x)_+]$. Assume the given sequent is reduced in the following:

\medskip  

$\enspace \enspace \displaystyle{\frac{\enspace \enspace \rightarrow G_1[\exists xA(x)_+, \exists xB(x)_+]} 
{\enspace \enspace \rightarrow G_1[\exists xA(x)_+, \exists xB(x)_+] \vee B(t) \enspace \enspace .}}$ 
$\enspace \enspace (\rightarrow \exists_+)$

\bigskip 

\noindent By Theorem 5.7, we obtain $\bf SN \it \vdash  \rightarrow \exists xA(x)$ or $\bf SN \it \vdash \rightarrow G_1[\; \;_+, \exists xB(x)_+] \vee B(t)$. 
If $\bf SN \it \vdash  \rightarrow \exists xA(x)$, then we find a term $s$ such that $\bf SN \it \vdash  \rightarrow A(s)$ by H.I., from 
which follows $\rightarrow G_1[A(s)_+, \exists xB(x)_+]$. Otherwise, there obtains 
$\bf SN \it \vdash \rightarrow G_1[\; \;_+, \exists xB(x)_+] \vee B(t)$, to which 
$\rightarrow G_1[\; \;_+, \exists xB(x)_+] \vee B(t)$ is reduced. Thus, $\bf SN \it \vdash \rightarrow G_1[A(s)_+, \exists xB(x)_+] \vee B(t)$ holds 
for any term $s$ by the thinning theorem.

Suppose that $\rightarrow G[\exists xA(x)_+]$ is reduced to $\rightarrow G[\exists xA(x)_+] \vee A(t)$. 
By Theorem 5.7, we obtain $\bf SN \it \vdash  \rightarrow \exists xA(x)$ or 
$\bf SN \it \vdash \rightarrow G[\; \;_+] \vee A(t)$. If $\rightarrow \exists xA(x)$ is provable, then 
we can find a term $t$ as required by H.I.. Otherwise, $\bf SN \it \vdash  \rightarrow G[A(t)_+]$ is obtained by the translation theorem, and 
this is the looked-for sequent. 

\it Proof of $\Leftarrow$ for Theorem 5.8\rm . Suppose that a term $t$ is found such that $\bf SN \it \vdash \enspace \rightarrow G[A(t)_+]$.
By the translation theorem, we have $\bf SN \it \vdash \enspace \rightarrow G[\; \;_+] \vee A(t)$. Then apply thinning theorem to it to have 
$\bf SN \it \vdash \enspace \rightarrow G[\exists xA(x)_+] \vee A(t)$, to which $\rightarrow G[\exists xA(x)_+]$ is reduced by 
$(\rightarrow \exists_+)$. $\Box$ 

Theorems 5.7 and 5.8 fully reflect the constructive feature of $\bf SN$. It is, of course, not satisfied in the classical theory. 

\begin{theorem} \rm (Equality theorems) \it 

\medskip 

$(1)$ \enspace \enspace $\bf SN \it \vdash F[s = t^+, A(s)^+] \rightarrow G[A(t)_+]$,

$(2)$ \enspace \enspace $\bf SN \it \vdash F[s = t^+, A(t)^+] \rightarrow G[A(s)_+]$,

$(3)$ \enspace \enspace $\bf SN \it \vdash F[s = t^+, A(s)^-] \rightarrow G[A(t)_-]$,

$(4)$ \enspace \enspace $\bf SN \it \vdash F[s = t^+, A(t)^-] \rightarrow G[A(s)_-]$,

$(5)$ \enspace \enspace $\bf SN \it \vdash F[s = t^+, A(s)^+, A(t)^-] \rightarrow \Delta$, 

$(6)$ \enspace \enspace $\bf SN \it \vdash F[s = t^+, A(t)^+, A(s)^-] \rightarrow \Delta$, 

\medskip 

\noindent where $s$ and $t$ are terms, and $A(s)$ and $A(t)$ are formulas obtained from $A(x)$  by substituing $s$ and 
$t$, respectively to $x$.
\end{theorem}

\it Proof\rm . (1)--(4) are proved simultanously. For the proof use is made of induction on the length of the formula $A(s)$.

\it Basis \rm Case 1: Suppose that $s = t$ are numerical and $s = t$ is true. In the case, $A(s)$ and $A(t)$ are formulas 
of length 0 equivalent to each other. Consequently, (1)--(4) are respectively provable in view of the axioms 6 and 7.

Case 2: Assume that $s = t$ is numerical and false. Then, (1)--(4) are all the instances of the axiom 2.

Case 3: Suppose that $s = t$ is not numerical. It is immediate that (1)--(4) respectively constitute the instances of the axiom 8. 

\it Induction steps \rm $A(s)$ is one of the forms $A_1(s) \wedge A_2(s)$, $A_1(s) \vee A_2(s)$, $A_1(s) \supset A_2(s)$, 
$\sim A_1(s)$, $\forall xA_1(x, s)$, $\exists xA_1(x, s)$. 

Let $A(s)$ be  $A_1(s) \wedge A_2(s)$. (1) is reduced as follows: 

\medskip  

\noindent $\displaystyle{\frac{F[s = t^+, A(s)_+] \rightarrow G[A_1(t) \wedge A_2(t)_+]}
{F[s = t^+, A(s)^+] \rightarrow G[A_1(t)_+]  | 
F[s = ^+t, A(s)^+] \rightarrow G[A_2(t)_+]  .}}$ 
$(\rightarrow \wedge_+)$

\bigskip 

\noindent The formulas under located are both provable in $\bf SN$ by H.I.. (1) is, hence, provable. 

Let $A(s)$ be  $A_1(s) \vee A_2(s)$. (1) is reduced as shown below:

\medskip  

\noindent $\displaystyle{\frac{F[s = t, A_1(s) \vee A_2(s)^+] \rightarrow G[A(t)_+]}
{F[s = t^+, A_1(s)^+] \rightarrow G[A(t)_+]  | 
F[s = t^+, A_2(s)^+] \rightarrow G[A(t)_+].}}$ 
$(\vee^+ \rightarrow)$

\bigskip 

\noindent By H.I., we finish it.

Let $A(s)$ be  $\sim A_1(s)$. So we have 
\medskip  

$F[s = t^+, \sim A_1(s)^+] \rightarrow G[\sim A_1(t)_+],$

\noindent which  is of the form: 

$F[s = t^+, A_1(s)^-] \rightarrow G[A_1(t)_-], \enspace \enspace  (*)$

\noindent and $(*)$ is provable by H.I.. Therefore, (1) is also provable.

Let $A(s)$ be $A_1(s) \supset A_2(s)$. It is suffices to take care of the following reduction:
$$\displaystyle{\displaystyle{
\displaystyle{\enspace F[s = t^+, A(s)^+] \rightarrow G[A_1(t) \supset A_2(t)_+] \enspace 
\over 
\displaystyle{\enspace F[s = t^+, A_1(s) \supset A_2(s))^+] \wedge A_1(t) \rightarrow A_2(t) \enspace } 
\enspace }\enspace (\rightarrow \supset_+) }
\over 
\displaystyle{\enspace F[s = t^+, A(s)^+] \wedge A_1(t) \rightarrow A_2(t) \vee A_1(s) | 
( F[s = t^+, A(s)^+] \wedge A_1(t)) \wedge A_2(s) \rightarrow A_2(t)  \enspace . }} \enspace (\supset^+ \rightarrow)$$
Hence, (1) is provable since the under located sequents are provable by H.I.. 

Let $A(s)$ be $\forall xA_1(x, s)$. (1) is reduced as shown below:

$$\displaystyle{\enspace F[s = t^+, \forall xA_1(x, s)^+] \rightarrow G[\forall xA_1(x, t)_+] \enspace 
\over 
\displaystyle{\enspace F[s = t^+, \forall A_1(x, s)^+] \rightarrow A_1(b, t) \enspace 
\over 
\displaystyle{\enspace F[s = t^+, A(s)^+] \wedge A_1(b, s) \rightarrow A_1(b, t) \enspace .
}  } \enspace  (\forall^+ \rightarrow) } \enspace (\rightarrow \forall_+)$$
It is obvious that the sequent at the bottom is provable by H.I.. (1), therefore, is provable. The remaining cases are proved similarly.  

We next give the proof of (5) and (6). It is remarked that (5) and (6) are proved concurrently like the proof of (1)--(4). (1)--(4) are also 
required for proving them. 

\it Basis \rm Case 1: Assume that $s = t$ is numerical and true. $A(s)$ and $A(t)$ are, then, equivalent formulas of length $0$. 
Consequently, (5) and (6) are, respectively, the instances of the axiom 5. 

Case 2: Suppose that $s = t$ is numerical and false. Then, (5) and (6) are the instances of the axiom 2 respectively.

Case 3: Assume that $s = t$ is not numerical. (5) and (6) are, then, the instances of the axiom 8.

\it Induction steps \rm We shall only prove some representative cases, while the remaining cases will ber analogously taken care of 
as in (1)--(4).

Let $A(s)$ be $A_1(s) \supset A_2(s)$. (5) is, then, reduced as follows: using $(\supset^+ \rightarrow)$, 
\medskip  

\noindent $\displaystyle{\frac{\enspace \enspace F[s = t^+, A_1(s) \supset A_2(s), A(t)^-] \rightarrow \Delta }
{F[s = t^+, A(s)^+, A(t)^-] \rightarrow \Delta \vee A_1(s) | 
F[s = t^+, A(s)^+, A(t)^-] \vee A_2(s) \rightarrow \Delta .}}$ 

\bigskip 

\noindent The left sequent reduced is provable Theorem 5.3.(2). The right sequent reduced is also  provable by H.I.. (5) is provable in 
$\bf SN$, thus.

Let $A(s)$ be $\sim A_a(s)$. (5), then, is thought of as the following sequent:

$F_1[s = t^+, A_1(s)^-, A_1(t)^+] \rightarrow \Delta  . \enspace \enspace (*)$

\noindent $(*)$ is provable  by H.I.. So is (5).

\begin{theorem} \mbox{} 

\medskip 

$(1)$ \enspace \enspace $\mbox{\bf SN} \vdash F[\; \;^-] \rightarrow G[A_+] \enspace \Rightarrow \enspace 
\mbox{\bf SN} \vdash F[A^-] \rightarrow G[\; \;_+] ,$

$(2)$ \enspace \enspace $\mbox{\bf SN} \vdash F[\; \;^+] \rightarrow G[A_-] \enspace \Rightarrow \enspace 
\mbox{\bf SN}  \vdash F[A^+] \rightarrow G[\; \;_-] , $
\end{theorem} 

\it Proof\rm . This is proved by induction on the length of the given tableau. $\Box$

Note that the converse does not hold in general. If it did, $\bf SN$ would be classical.

\begin{theorem} \mbox{}

\medskip

$(1)$ \enspace \enspace $\mbox{\bf PCN} \vdash F[A^+] \rightarrow \enspace \Rightarrow 
\enspace \mbox{\bf SN} \vdash F[\; \;^+] \rightarrow \enspace \sim A .$

$(2)$ \enspace \enspace $\mbox{\bf PCN} \vdash F[A^-] \rightarrow \enspace \Rightarrow \enspace \mbox{\bf SN} \vdash F[\; \;^+] \rightarrow A .$

\end{theorem} 

\it Proof\rm . This is proved by induction on the length of the given tableau. $\Box$

\begin{theorem} \mbox{} 

\medskip 

$(1)$ \enspace \enspace $\mbox{\bf FN} \vdash F[P^+] \rightarrow \enspace \Rightarrow \enspace 
\mbox{\bf FN} \vdash F[\; \;^+] \rightarrow \enspace \sim P ,$

$(2)$ \enspace \enspace $\mbox{\bf FN}  \vdash F[P^-] \rightarrow \enspace \Rightarrow \enspace 
\mbox{\bf FN} \vdash F[\; \;^-] \rightarrow P, $

\medskip 

\noindent where $P$ is a prime formula.
\end{theorem} 

\it Proof\rm . This is proved by induction on the length of the given tableau. $\Box$

We obtain the following important theorem, again, after Sch\"utte [1977] \cite{schutte1977}.

\begin{theorem} \rm (Extended axiom theorems) \it 

\medskip 

$(1)$ \enspace \enspace $\bf SN \it \vdash F[A^+] \rightarrow G[B_+]$,

$(2)$ \enspace \enspace $\bf SN \it \vdash F[A^-] \rightarrow G[B_-]$,

$(3)$ \enspace \enspace $\bf SN \it \vdash F[A^+, B^-] \rightarrow \Delta$,

\medskip 

\noindent where $A$ and $B$ are equivalent formulas.
\end{theorem}

\it Proof\rm . (1) and (2) are proved simultaniously before we proceed to the proof of (3). By induction on the length of 
formula $A$, (1) and (2) are proved like the proof of Theorem 5.3. (3) is proved similarly with the help of (1) and (2). 
Note that Theorem 5.13.(1)--(2) could be proved on the basis of Theorem 5.3 and the $\bf SN$-cut elimination to be proved in \S 9 
as an analogue of the corresponding one of Sch\"utte's $\Delta_1^1$-analysis. (For the detais, refer to Sch\"utte [1977] \cite{schutte1977}.) 

Lastly, we are presenting the following theorems, which though, trivial, will be made use of in \S 5. 

\begin{theorem} \rm (Reduction theorem) \it 

\medskip 

\enspace \enspace $\bf FN \it \vdash \Gamma \rightarrow \enspace \Rightarrow \enspace \bf PCN \it \vdash \Gamma\rightarrow \enspace$.

\medskip 
\end{theorem}

\it Proof\rm . Trivial. $\Box$

\begin{theorem} \rm \mbox {} \it 

\medskip 

\enspace \enspace $\bf FN \it \vdash F[A \supset B^+] \rightarrow \enspace \Rightarrow \enspace \bf PCN \it \vdash F[\; \;^+] \rightarrow \enspace $.

\medskip 

\end{theorem}

\it Proof\rm . The theorem is easily proved by induction on the length of the given tableau. $\Box$ 

\begin{lemma} \mbox{} 

\medskip 

$(1)$ \enspace \enspace $\mbox{\bf SN} \vdash F[A^+] \rightarrow \Delta \enspace \Leftrightarrow \enspace 
\mbox{\bf SN} \vdash F[\; \;^+] \wedge A \rightarrow \Delta ,$

$(2)$ \enspace \enspace $\mbox{\bf SN} \vdash F[A^-] \rightarrow \Delta \enspace \Leftrightarrow \enspace 
\mbox{\bf SN} \vdash F[\; \;^-] \wedge \sim A \rightarrow \Delta ,$

$(3)$ \enspace \enspace $\mbox{\bf SN} \vdash \Gamma \rightarrow G[A_+]  \enspace \Leftrightarrow \enspace 
\mbox{\bf SN} \vdash \Gamma \rightarrow G[\; \;_+] \vee A ,$

$(4)$ \enspace \enspace $\mbox{\bf SN} \vdash \Gamma \rightarrow G[A_-]  \enspace \Leftrightarrow \enspace 
\mbox{\bf SN} \vdash \Gamma \rightarrow G[\; \;_-] \vee \sim A ,$
\end{lemma} 

\it Proof\rm . The lemma is easily proved by induction on the length of the given tableau. Or just apply Theorem 5.5. If you prefer a semantical proof, 
then take simultaneous induction on the numbers of prodedures which determins the positive and negative parts. $\Box$


\section{Two cut elimination theorems of $\bf SN$}

There are two cut elimination theorems for the proposed $\bf SN$. One is a natural orthodox cut elimination theorem. The other is author's 
new type cut elimination theorem. First we are presenting the former.  

\begin{theorem} \rm ($\bf SN$-cut elimination theorem)  (For the simpler predicate logic case, refer to Kanai [1984] \cite{kanai}) \it 
$$(\mbox{\bf SN}  \vdash \Gamma \rightarrow G[A_\pm]  \mbox{ and  }  \mbox{\bf SN} \vdash F[A^\pm] \rightarrow \Delta) 
\enspace \Rightarrow \enspace \mbox{\bf SN} \vdash \Gamma \wedge F[\; \;^\pm] \rightarrow G[\; \;_\pm] \vee \Delta \enspace.$$
\end{theorem}

The proof of Theorem 6.1 will be presented in \S 9. 

For proving the consistency of $\bf SN$, we need the $\bf SN$-cut elimination theorem. Assuming 
$\bf SN \it \vdash \enspace \rightarrow A$ and \enspace 
$\bf SN \it \vdash \enspace \rightarrow \enspace\sim A$ for some formula $A$ of $\bf SN$, we obtain $\bf \it \vdash \rightarrow$ by the 
application of the $\bf SN$-cut elimination theorem. The sequent $\rightarrow$, however, is not provable in $\bf SN$, since we can not 
reduce it to any other sequent. This leads us to a contradiction. In fact, if $\bf SN \it \vdash \rightarrow$, we have 
$\bf SN \it \vdash \enspace \rightarrow A$ and \enspace 
$\bf SN \it \vdash \enspace \rightarrow \enspace \sim A$ for some formula $A$ of $\bf SN$ by the thinning theorem. The consistency of 
$\bf SN$ could be understood more deeply on the basis of the above discussion. In fact, $\bf SN$ is obviously consistent, since $\bf SN$ is 
formulated as a cut-free tableau system. The sequent $\rightarrow 1 = 0$ is not provable in $\bf SN$.


\begin{theorem} \rm (Consistency theorem of $\bf SN$) \it 

\medskip 

\enspace \enspace $\bf SN$ is  consistent. 

\end{theorem}

\it Proof\rm . Suppose, if possible, $\bf SN$ be inconsistent. Thus,  $\bf SN \it \vdash \enspace \rightarrow A$ and \enspace 
$\bf SN \it \vdash \enspace \rightarrow \enspace\sim A$ for some formula $A$ of $\bf SN$.  By the above discussion, 
$\rightarrow$ is provable in $\bf SN$. This is impossible. $\Box$

\begin{theorem} The following sentences are equivalent to each other. 

\medskip 

$(1)$ \enspace \enspace $\bf SN$ is inconsistent.

$(2)$ \enspace \enspace $\bf SN \it \vdash \enspace \rightarrow A$ and \enspace $\bf SN \it \vdash \enspace \rightarrow \enspace\sim A$ 
for some formula $A$ of $\bf SN$.

$(3)$ \enspace \enspace $\bf SN \it \vdash \enspace \rightarrow \enspace $.

$(4)$ \enspace \enspace $\bf SN \it \vdash \enspace \rightarrow 1 = 0$, where $1$ is the abbreviation of $0'$.

$(5)$ \enspace \enspace $\bf SN \it \vdash \enspace \rightarrow A$ for any formula $A$ of $\bf SN$.

\end{theorem}

\it Proof\rm . (1) $\Leftrightarrow$ (2) $\Leftrightarrow$ (3) is clear from Theorem 6.2 and the discussion above 
theorem 6.2 of (3) $\Rightarrow$ (4), (4) is obtaind from (3) by Thinning theorem. Conversely, assume (4). 
Sequent $1 = 0 \rightarrow $ is provable in $\bf SN$, since it is the axiom 2. From (4) and $1 = 0 \rightarrow$, 
we obtain by the $\bf SN$-cut elimination theorem. Of (5) $\Rightarrow$ (4), (4) is easily obtained from (5). We 
are taking care of (3) $\Rightarrow$ (5), lastly. For any formula $A$, $\rightarrow A$ is provable, since it is 
obtained from $\enspace \rightarrow \enspace$ by Thinning theorem. $\Box$ 

Now, another cut elimination theorem will be presented.

\begin{theorem} \rm ($\bf PCN$-cut elimination theorem, Inou\'e [1984] \cite{inouen}) \it 
$$(\mbox{\bf PCN} \vdash F_1[A^+] \rightarrow \mbox{and }  \mbox{\bf PCN} \vdash F_2[A^-] \rightarrow \enspace) \Rightarrow \enspace 
\mbox{\bf PCN} \vdash F_1[\; \;^+] \wedge F_2[\;\;^-] \rightarrow \enspace .$$
\end{theorem}

It is again emphasized that the  $\bf PCN$-cut elimination theorem plays an important part of our embedding of $\bf CN$ 
(classical number theory without complete induction) in $\bf SN$ to be proved in \S 8. The theorem successfully takes care 
of the modus ponens of $\bf CN$ within the bound of $\bf PCN$, since it is a version of the modus ponens in $\bf PCN$. 
This is understandable if we remember that $\bf PCN$ is just like the classical system with respect to its axioms and 
reduction rules. In addition, it is noted, the $\bf PCN$-cut elimination theorem is by no means generalized to $\bf SN$ as 
easily shown by a counterexample. (For the details, refre to Inou\'e [1984] \cite{inouen}. It was conjectured in the paper 
that the new cut elimination theorem could be proved in the constructive predicate logic involving strong negation without 
any restriction.) 

\it Proof of Theorem 6.4\rm . For proving this we use a conventional double induction on the grade and rank of the 
$\bf PCN$-cut (i.e., double induction on $\omega \cdot g + r$ with $g, r < \omega$). 

The formula $A$ occurring in the premises of the theorem is called \it cut formula\rm . The \it grade \rm of a cut 
formula $A$ denoted by $g(A)$, is the length of the formula $A$. 

The \it left rank \rm of the cut with the cut formula $A$, denoted by $rank_l(A)$, is the sum of sequents which contain $A$, 
beginning with the left premise up to that in which $A$ is introduced first. The \it right rank \rm of the cut with the cut formula 
$A$, denoted by $rank_r(A)$ is similarly defined. The sum of $rank_l(A)$ and 
$rank_r(A)$ is called the rank of the cut with the cut formula $A$, and denoted by $rank(A)$. Notice that $g(A) \geq 0$ and 
$rank(A) \geq 0$. 

For brevity we shall use the following abbreviations, namely, $Tr$, $NTr$, $Fal$ and $Eq$ which, respectively, stand for:

\smallskip 

$Tr(A)$ $\Leftrightarrow$ ($A$ is a true constant prime formula),

$NTr(A)$ $\Leftrightarrow$ ($A$ is not a ture constant prime formula),

$Fal(A)$ $\Leftrightarrow$ ($A$ is a false constant prime formula),

$Eq(A, B)$ $\Leftrightarrow$ ($A$ and $B$ are equivalent formulas of length $0$).

\smallskip 

By Theorem 5.14, it suffices to consider only four axioms (1, 2, 5 and 8) and six reduction rules ($(\wedge^- \rightarrow)$, 
$(\vee^+ \rightarrow)$, $(\forall^+ \rightarrow)$,  $(\forall^- \rightarrow)$, 
$(\exists^+ \rightarrow)$ and $(\exists^- \rightarrow)$) for the construction of the tableau for the premises. 

\medskip 

Case 1 ($g(A) = 0$ and $rank(A) = 2$): 

Subcase A: Let $F_1$ be $F'_1[A^+, B^-]$. Assume $Eq(A, B)$ and $Tr(A)$. $F'_1[A^+, B^-] \wedge 
F'_1[\; \;^+, B^-] \wedge F_2[\; \;^-] \rightarrow$ is the axiom 1 since $Tr(B)$ by  the assumption.

Subcase B:  Let $F_1$ be $F'_1[A^+, B^-]$ and  $F_2$ be $F'_2[A^-, C^+]$.  Assume $Eq(A, B)$, $Eq(A, C)$ and $NTr(A)$. 
$F'_1[\; \;^+, B^-] \wedge F'_2[\; \;^-, C^+] \rightarrow$ is, then, the axiom 6 since $Eq(B, C)$ by the assumption. 

Subcase C:  Let $F_1$ be $F'_1[A^+, B^-]$ and  $F_2$ be $F'_2[A^-, C^+]$.  Assume $Eq(A, B)$, $Eq(A, C)$ and $Tr(C)$. 
$F'_1[\; \;^+, B^-] \wedge F'_2[\; \;^-, C^+] \rightarrow$ is, then, the axiom 1 by $Tr(C)$. 

Subcase D:  Let $F_1$ be $F'_1[A^+, B^-]$ and  $F_2$ be $F'_2[A^-, C^+]$.  Assume $Eq(A, B)$, $NTr(A)$ and $Fal(C)$. 
$F'_1[\; \;^+, B^-] \wedge F'_2[\; \;^-, C^+, D^-] \rightarrow$ is, then, the axiom 2 by $Fal(C)$. 

Subcase E:  Let $F_1$ be $F'_1[A^+, B^-]$ and  $F_2$ be $F'_2[A^-, C^+, D^-]$.  Assume $Eq(A, B)$ and $Eq(C, D)$. 
$F'_1[\; \;^+, B^-] \wedge F'_2[\; \;^-, C^+] \rightarrow$ is, then, the axiom 5 by $Eq(C, D)$.

Subcase F:  Let $F_1$ be $F'_1[A^+, B^-]$ and  $F_2$ be $A(x_1, \dots ,x_n)$ $(n \geq 1)$. Assume $Eq(A, B)$ and 
$A(x_1, \dots ,x_n)$ satisfies the condition of the axiom 8. The case has already been taken of 
by the subcases A--E of the case 1.

Subcase G: Let $F_2$ be $F'_2[A^-, B^+]$. Assume $Fal(A)$ and $Eq(A, B)$. $F_1[\; \;^+] \wedge 
F'_2[\; \;^-, B^+] \rightarrow$ is the axiom 2 since we have $Fal(B)$ from the assumption. 

Subcase H: Let $F_2$ be $F'_2[A^-, B^-]$. Assume $Fal(A)$ and $Tr(B)$. $F_1[\; \;^+] \wedge 
F'_2[\; \;^-, B^-] \rightarrow$ is the axiom 1 since we have $Tr(B)$ from the assumption. 

Subcase I: Let $F_2$ be $F'_2[A^-, B^+]$. Assume $Fal(A)$ and $Fal(B)$. $F_1[\; \;^+] \wedge 
F'_2[\; \;^-, B^+] \rightarrow$ is the axiom 2 by $Fal(B)$. 

Subcase J:  Let $F_2$ be $F'_2[A^-, B^+, C^-]$.  Assume $Fal(A)$ and $Eq(B, C)$. 
$F_1[\; \;^+] \wedge F'_2[\; \;^-, B^+, C^-] \rightarrow$ is, then, the axiom 5 by $Eq(B, C)$. 

Subcase K:  Let $F_2$ be $A(x_1, \dots ,x_n)$ $(n \geq 1)$.  $A(x_1, \dots ,x_n)$ satisfies the condition of the axiom 8. 
The case has already been considered by the subcases G--K of the case 1. 

The case that $F_1[A^+] \rightarrow$ is the the axiom 2 by $Fal(S)$ has been already taken care of by the subcases 
G--K of the case1. 

Subcase L: Assume that the cut formula $A$ of $F_1$ is not the principal part of an axiom. From the assumption, 
 $F_1[\; \;^+] \rightarrow$ is one of the axioms 1, 2, 5 and 8.  $F_1[\; \;^+] \wedge F_2[\; \;^-] \rightarrow$ is 
also an axiom. 

Subcase M:  Let $F_1$ be $A(x_1, \dots ,x_n)$ $(n \geq 1)$.  $A(x_1, \dots ,x_n)$ satisfies the condition of the axiom 8. 
The case has already been considered by the above subcases A--L of the case 1. 

Case 2 ($g(A) = 0$, $rank_l(A) = 1$ and $rank_r(A) > 1$): By $g(A) = 0$, the cut formula of the premises is not the 
principal formula of any reduction rules. There are six reduction rules to be taken up, but only two typical ones 
will be taken care of. 

Let $F_2$ be $F'_2[A^-, B \wedge C^-]$. If $F_2$ is reduced by $(\wedge^- \rightarrow)$,  we, then, obtain two sequents 
$F'_2[A^-, B^-] \rightarrow $ and $F'_2[A^-, C^-] \rightarrow $. By H.I., 

$\bf PCN \it \vdash F_1[\; \;^+] \wedge F'_2[\; \;^-, B^-] \rightarrow$, \enspace (2--1)

$\bf PCN \it \vdash F_1[\; \;^+] \wedge F'_2[\; \;^-, C^-] \rightarrow$, \enspace (2--2)

\noindent By (2--1) and (2--2), $F_1[\; \;^+] \wedge F'_2[\; \;^-, B \wedge C^-] \rightarrow$ is provable in $\bf PCN$, 
since it is demonstrable by the following reduction:

\medskip  

$\displaystyle{\frac{F_1[\; \;^+] \wedge F'_2[\; \;^-, B\wedge C^-] \rightarrow]}
{F_1[\; \;^+] \wedge F'_2[\; \;^-, B^-] \rightarrow \enspace  | \enspace
F_1[\; \;^+] \wedge F'_2[\; \;^-, C^-] \rightarrow  .}}$ 
$(\wedge^- \rightarrow)$ 

\bigskip 

\noindent Another case we wish to try is the one where $F_2$ is $F'_2[A^-, \forall xB(x)^-]$. By the inversion 
theorem $F_2$ is reduced, and we obtain the sequent $F'_2[\; \;^-, B(b)^-] \rightarrow$, where $b$ is a proper 
(eigen) variable. By H.I., 
$$\bf PCN \it \vdash F_1[\; \;^+] \wedge F'_2[\; \;^-, B(b)^-] \rightarrow , \enspace \mbox{(2--3)}$$
Let $c$ be a variable not occurring as a free variagle in $F_1[\; \;^+] \wedge F'_2[\; \;^-,  \forall xB(x)^-] \rightarrow$. 
Let $x_2, \dots , x_n$ are pairwise distinct free variables of $F_1[\; \;^+] \wedge F'_2[\; \;^-,  B(b)^-] \rightarrow$. 
Assume further that $x_2, \dots , x_n$ are not bound and distinct from $b$.  By Theorem 5.1 and (2--3), 
$$\bf PCN \it \vdash F_1[\; \;^+] \wedge F'_2[\; \;^-, B(b)^-] [c/b, x_2/x_2, \dots , x_n/x_n] \rightarrow  .$$
Namely, we have 
$$\bf PCN \it \vdash F_1[\; \;^+] \wedge F'_2[\; \;^-, B(b)^-]  \rightarrow . \enspace \mbox{(2--5)}$$
From (2--5), we obtain $\bf PCN \it \vdash F_1[\; \;^+] \wedge F'_2[\; \;^-, \forall xB(x)^-]  \rightarrow$, since it is 
reduced as follows: 

\medskip  

$\displaystyle{\frac{F_1[\; \;^+] \wedge F'_2[\; \;^-, \forall xB(x)^-]  \rightarrow}
{F_1[\; \;^+] \wedge F'_2[\; \;^-, B(b)^-] \rightarrow  .}}$ 
$(\forall^- \rightarrow)$ 

\bigskip 

\noindent Notice that we employ Theorem 5.6 (Interchange theorems) in the case of $(\forall^+ \rightarrow)$ and 
$(\exists^- \rightarrow)$, 

Case 3 ($g(A) = 0$, $rank_l(A) > 1$ and $rank_r(A) = 1$): The proof is similar to that for the case 2. 

Case 4 ($g(A) = 0$, $rank_l(A) > 1$ and $rank_r(A) > 1$): The treatment of this case is analogous to that for the case 2. 

Case 5 ($g(A) > 0$ and  $rank(A) = 2$): 

Subcase A: Suppose that the cut formula $A$ is of the form $A_1 \wedge A_2$. By the inversion theorem, 
$F_2[A_1^-] \rightarrow$ and $F_2[A_2^-] \rightarrow$ are obtained from $F_2[A_1 \wedge A_2^-] \rightarrow$. In what 
follows, such inferences as: 
$$(\vdash \Gamma_1 \rightarrow \Delta_1 \mbox { and } \vdash \Gamma_2 \rightarrow \Delta_2) \Rightarrow \enspace
\vdash \Gamma_3 \rightarrow \Delta_3$$
is described \'a la Gentzen as:
$$\displaystyle{\frac{\enspace \enspace \Gamma_1 \rightarrow \Delta_1  \enspace \enspace \enspace \enspace 
\Gamma_2 \rightarrow \Delta_2 \enspace \enspace} 
{\Gamma_3 \rightarrow \Delta_3 \enspace  ,}}$$
where  $\Gamma_1 \rightarrow \Delta_1$ and $\Gamma_2 \rightarrow \Delta_2$ are understood accompanying the 
respective closed tableau to be obtained by reducing them. (A \it closed \rm tableau is one of which each 
branch ends with one of axioms. A sequent which has at least one closed tableau is provable.) We, then, obtain a 
series of inferences as follows:
$$\displaystyle{\enspace \mbox{Interchange } 
\displaystyle{ \mbox{H.I. on $grade$} \enspace \displaystyle{ 
\displaystyle{\enspace F_1[A_1 \wedge A_2^+] \rightarrow \enspace \enspace  F_2[A_1^-] \rightarrow  \enspace 
\over 
\displaystyle{\enspace (F_1[\; \;^+] \wedge A_2) \wedge F_2[\; \;^-] \rightarrow \enspace } 
\enspace }\enspace }
\over 
\displaystyle{\enspace (F_1[\; \;^+] \wedge F_2[\; \;^-]) \wedge A_2 \rightarrow 
 }} \enspace \enspace 
\enspace \enspace F_2[A_2^-] \rightarrow \enspace 
\over 
\displaystyle{\enspace (F_1[\; \;^+] \wedge F_2[\; \;^-]) \wedge F_2[\; \;^-] \rightarrow  \enspace 
\over 
\displaystyle{\enspace F_1[\; \;^+] \wedge F_2[\; \;^-] \rightarrow \enspace .
}  } \enspace  \mbox{Contraction} } \enspace \mbox{H.I. on $grade$}$$
The interchange theorem is called for with a view to making cut more easily applicable. 

Subcase B: Assume that the cut formula has the form $A_1 \vee A_2$.  By the inversion theorem, we obtain 
$F_1[A_1^+] \rightarrow$ and $F_1[A_2^+] \rightarrow$ from $F_1[A_1 \vee A_2^+] \rightarrow$. We, then, have: 
$$\displaystyle{\enspace F_1[A_2^+] \rightarrow  \enspace \mbox{Translation (twise) } 
\displaystyle{ \mbox{H.I. on $grade$} \enspace \displaystyle{ 
\displaystyle{\enspace  F_1[A_1^+] \rightarrow  \enspace \enspace F_2[A_1 \vee A_2^-] \rightarrow \enspace 
\over 
\displaystyle{\enspace (F_1[\; \;^+]\wedge F_2[\; \;^-]) \wedge \sim A_2 \rightarrow \enspace } 
\enspace }\enspace }
\over 
\displaystyle{\enspace (F_1[\; \;^+] \wedge F_2[\; \;^-]) \wedge A_2 \rightarrow 
 }}  
\over 
\displaystyle{\enspace F_1[\; \;^+] \wedge (F_1[\; \;^+]) \wedge F_2[\; \;^-]) \rightarrow  \enspace 
\over 
\displaystyle{\enspace F_1[\; \;^+] \wedge F_2[\; \;^-] \rightarrow \enspace .
}  } \enspace  \mbox{Contraction} } \enspace \mbox{H.I. on $grade$}$$

Subcase C: Suppose that the cut formula $A$ is $A_1 \supset A_2$. By Theorems 5.13 and 5.14, 
$\bf PCN \it \vdash F_1[\; \;^+] \rightarrow$ is forthcoming from $F_1[\; \;^+]$ is forthcoming from 
$F_1[A_1\supset A_2^+] \rightarrow$. Then $\bf PCN \it \vdash F_1[\; \;^+] \wedge F_2[\; \;^-] \rightarrow$ is 
obtained by the thinning theorem. 

Subcase D: Let $A$ be $\sim A$. $F_1[A_1^+]$ and $F_2[A_1^-]$ are thought of as  $F_3[A_1^-]$ and $F_4[A_1^+]$, 
respectively. We, then, obtain the following series of inferences:

$$\displaystyle{\enspace F_4[A_1^+] \rightarrow \enspace F_3[A_1^-] \rightarrow 
\over 
\displaystyle{\enspace F_4[\; \;^+] \wedge F_3[\; \;^-] \rightarrow \enspace 
\over 
\displaystyle{\enspace F_3[\; \;^+] \wedge F_4[\; \;^-] \rightarrow  \enspace .
}  } \enspace  \mbox{Interchange} } \enspace \mbox{H.I. on $grade$}$$
$F_3[\; \;^+] \wedge F_4[\; \;^-]$ is obviously $ F_1[\; \;^+] \wedge F_2[\; \;^-]$.

Subcase E: Assume that the cut formula $A$ is of the form $\forall xA(x)$. $F_1[\forall xA(x)^+] \rightarrow$ is an axiom 
by the assumption that $rank_l(\forall xA(x)) = 1$.  $F_1[\; \;^+] \rightarrow$ is also ana axiom since $\forall xA(x)$ is 
by no means the principal part of an axiom. $F_1[\; \;^+] \wedge F_2[\; \;^-] \rightarrow$ remains an axiom, thus.

Subcase F: Let the cut formula $A$ be $\exists xA(x)$. $F_2[\; \;^-] \rightarrow$ is an axiom again by the assumption. It 
follows that $F_1[\; \;^+] \wedge F_2[\; \;^-] \rightarrow$ remains an axiom of the same kind.

Case 6 ($g(A) = 0$, $rank_l(A) = 1$ and $rank_r(A) > 1$): 

Subcase A: Assume $A$ is of th form $A_1 \wedge A_2$. By the assumption, $F_2[A_1 \wedge A_2^-] \rightarrow $ is not an axiom. 
If $F_2[A_1 \wedge A_2^-] \rightarrow $ is reduced by ($\wedge^- \rightarrow$), we have: 

\medskip  

$\displaystyle{\enspace \enspace \frac{F'_2[A_1 \wedge A_2^-, B \wedge C^-] \rightarrow \enspace \enspace }
{\enspace \enspace F'_2[A_1 \wedge A_2^-, B^-] \rightarrow \enspace | \enspace F'_2[A_1 \wedge A_2^-, C^-] \rightarrow \enspace \enspace ,}}$ 
$(\wedge^- \rightarrow)$ 

\bigskip 

\noindent $F'_2$ is $F_2$. We, then, have: 

\medskip  

$\displaystyle{\enspace \enspace \frac{\enspace \enspace F_1[A_1 \wedge A_2^+] \rightarrow \enspace \enspace F'_2[A_1 \wedge A_2^-, B^-] 
\rightarrow \enspace \enspace }
{\enspace \enspace F_1[\; \;^+] \wedge F'_2[\; \;^-, B^-] \rightarrow  \enspace \enspace , \enspace \enspace (*)}}$ 
$\enspace \mbox{H.I. on $rank$}$ 
 
\bigskip 

$\displaystyle{\enspace \enspace \frac{\enspace \enspace F_1[A_1 \wedge A_2^+] \rightarrow \enspace \enspace F'_2[A_1 \wedge A_2^-, C^-] 
\rightarrow \enspace \enspace }
{\enspace \enspace F_1[\; \;^+] \wedge F'_2[\; \;^-, C^-] \rightarrow  \enspace \enspace . \enspace \enspace (**)}}$ 
$\enspace \mbox{H.I. on $rank$}$ 

\bigskip 

\noindent $F_1[\; \;^+] \wedge F'_2[\; \;^-, B \wedge C^-] \rightarrow$ is, thus, reduced to (*) and (**). The case that $(\vee^+ \rightarrow)$, 
$(\forall^+ \rightarrow)$ and $(\forall^- \rightarrow)$ are employed for the reduction are similarly taken care of as in the above case. 

Subcase B: Let $A$ be $A_1 \vee A_2$. By $g(A) > 1$, $F_1[\;\;^+] \rightarrow$ is an axiom. Then, $F_1[\; \;^+] \wedge F_2[\; \;^-] \rightarrow$ 
also constitutes an axiom.

Subcase C: Assume the cut formula $A$ is $A_1 \supset A_2$. The case is analogously dealt with as in the subcase $C$ of the case 5.

Subcase D: Assume that the cut formula $A$ is one of forms $\forall xA(x)$ and $\exists xA(x)$. In each case, we similarly obtain 
$\bf PCN \it \vdash F_1[\; \;^+] \wedge F_2[\; \;^-] \rightarrow$ as in the preceding subcase $B$ of the case 6. 

Subcase E: let $A$ be $\sim A_1$. The case is analogously proved as in the subcase D of the case 4. 

Case 7 ($g(A) > 1$, $rank_l(A) > 1$ and $rank_r(A) = 1$): The case is analogously treated as in the case 6. 

Case 8 ($g(A) > 1$, $rank_l(A) > 1$ and $rank_r(A) > 1$): 

Subcase A: If the cut formula $A$ is one of $A_1 \wedge A_2$, $A_1 \vee A_2$, $A_1 \supset A_2$ and $\sim A_1$. Then each case is 
similarly taken care of as in the case 5. 

Subcase B: Let $A$ be $\forall xA_1(x)$. Assume that $F_1[\forall xA_1(x)^+] \rightarrow$ is reduced in the following way: 

\medskip  

$\displaystyle{\enspace \enspace \frac{F_1[\forall xA(x)^+] \rightarrow \enspace \enspace }
{\enspace \enspace F_1[\forall xA(x)^+] \wedge A_1(t) \rightarrow  \enspace \enspace .}}$ 
$(\wedge^- \rightarrow)$ 

\bigskip 

\noindent We, then, have: 
$$\displaystyle{\enspace F_1[\forall xA_1(x)^+] \wedge A_1(t) \rightarrow \enspace F_2[\forall xA_1(x)^-] \rightarrow 
\over 
\displaystyle{\enspace (F_1[\; \;^+] \wedge A_1(t)) \wedge F_2[\; \;^-] \rightarrow \enspace 
\over 
\displaystyle{\enspace (F_1[\; \;^+] \wedge F_2[\; \;^-]) \wedge A_1(t) \rightarrow  \enspace .
}  } \enspace  \mbox{Interchange} } \enspace \mbox{H.I. on $rank$}$$
Now, by the inversion theorem there obtains $F_2[A_1(b)^-] \rightarrow$ from $F_2[\forall xA_1(x)^-] \rightarrow$. 
From this we obtain by substitution $F_2[A_1(t)^-] \rightarrow$. We, then, have: 

$$\displaystyle{\enspace (F_1[\; \;^+] \wedge F_2[\; \;^-]) \wedge A_1(t) \rightarrow \enspace F_2[A_1(t)^-] \rightarrow 
\over 
\displaystyle{\enspace (F_1[\; \;^+] \wedge F_2[\; \;^-]) \wedge F_2[\; \;^-] \rightarrow \enspace 
\over 
\displaystyle{\enspace (F_1[\; \;^+] \wedge F_2[\; \;^-]) \wedge A_1(t) \rightarrow  \enspace .
}  } \enspace  \mbox{Contraction} } \enspace \mbox{H.I. on $grade$}$$
If the cut formula $A$ is not principal formula of reduction rules, them we prove the cases by H.I. on rank as in the 
subcase A of the case 6.

Subcase C: Let $A$ be $\exists xA_1(x)$. The treatment is similar to that for the subcase B of the case 8. $\Box$


\section{Embedding theorem of $\bf CN$ in $\bf SN$}

In this section, we will proe a Kolmogorov-G\"odel-type embedding theorem of $\bf CN$ in $\bf SN$, which is the core for the 
consistency proof in \S 8. The following lemmas will be presented first for facilitating the proof. 

\begin{lemma} If a sequent $\Gamma \rightarrow \enspace$ is provable in $\bf SN$ and contains no implication sign, then 
$\Gamma \rightarrow \enspace$ is provable in $\bf PCN$.
\end{lemma}

\noindent \it Proof\rm . $\Gamma \rightarrow$ is provable in $\bf FN$ by the premise. The sequent is, then, provable in $\bf PCN$  by 
Theorem 5.14. $\Box$

\begin{theorem} \mbox{}

\medskip 

$(1)$ \enspace \enspace $\bf SN \it \vdash A \rightarrow B \enspace \Leftrightarrow \enspace 
\bf SNH \rm \vdash A \supset B$,

$(2)$ \enspace \enspace $\bf SN \it \vdash A \rightarrow \enspace \Leftrightarrow \enspace 
\bf SNH \rm \vdash A \supset \enspace \sim A$,

$(3)$ \enspace \enspace $\bf SN \it \vdash \enspace \rightarrow A \enspace \Leftrightarrow \enspace 
\bf SNH \rm \vdash A$,

$(4)$ \enspace \enspace $\bf SN \it \vdash \enspace \rightarrow \enspace \Leftrightarrow \enspace 
\bf SNH \rm \vdash A$ for every formula $A$ of $\bf SN$,

\medskip 

\noindent where $A$ and $B$ are formulas of $\bf SN$ $($\mbox{or} $\bf SNH$ $)$, and $\bf SNH$ is a Hilbert-type version of $\bf SN$.

\end{theorem}

\noindent \it Proof\rm . An easy proof will be omitted. Although the theorem is made use of in the proof of $\Leftarrow$ in the 
following embedding theorem, $\Leftarrow$ in the theorem  is dispensed with for the consistency proof of $\bf CN$. (cf. 
Ishimoto[1970] \cite{ishimotos} amd Takano [1970] \cite{takano}) $\Box$

\begin{theorem} \rm (Embedding theorem, for a simpler predicate logic part, Ishimoto [197?] \cite{ishimotoe} 
(refer to \cite{inoueac}, \cite{akama})) \it 
$$\mbox{\bf CN} \vdash^n A \enspace \Leftrightarrow \enspace 
\mbox{\bf PCN}\vdash \enspace \sim TA \rightarrow \enspace .$$where $A$ is a formula of $\bf CN$ and the operator $T$, 
which translates a formula of $\bf CN$ into its counterpart  in $\bf SN$, is recursively defined as follows:

\medskip 

$TA = A$ for every prime formula $A$ of $\bf CN$,

$T\sim A = \enspace \sim TA$,

$TA \wedge B = TA \wedge TB$,

$TA \vee B = TA \vee TB$,

$TA\supset B = \enspace \sim TA \vee TB$,

$T \forall xA(x) = \forall xTA(x)$,

$T \exists xA(x) = \exists xTA(x)$.

\end{theorem}

Let us consider the following example of the translation:

\smallskip 

\noindent $ \enspace \enspace T\forall y \sim \exists x((x = 0 \wedge y = 0') \supset x = y)$, 

\noindent $= \forall y T\sim \exists x((x = 0 \wedge y = 0') \supset x = y)$, 

\noindent $= \forall y \sim T\exists x((x = 0 \wedge y = 0') \supset x = y)$, 

\noindent $= \forall y \sim \exists xT((x = 0 \wedge y = 0') \supset x = y)$, 

\noindent $= \forall y \sim \exists x(\sim T(x = 0 \wedge y = 0') \vee Tx = y)$, 

\noindent $= \forall y \sim \exists x(\sim (Tx = 0 \wedge Ty = 0') \vee x = y)$, 

\noindent $= \forall y \sim \exists x(\sim (x = 0 \wedge y = 0') \vee x = y)$.

\smallskip 

\noindent It is noticed here that the logical symbols in the formula to be translated 
are different from the original ones which are, of course, interpreted classically. They are 
constructive as seeen from Theorems 5.7 and 5.8. It is remarked $TA$ \it has no implication sign \rm as 
seen from the definition of $T$ above.

\noindent \it Proof of Theorem 7.3\rm . We are first taking care of $\Rightarrow$ of the theorem 
by induction on the length $n$ of the proof of $A$. 

\it Basis \rm $(n = 0)$ 

Case 1 $(\mbox{\bf CN} \vdash^0 A \supset . B \supset A)$: The formula is translated by $T$ into the following 
sequent of $\bf SN$:
$$\sim (\sim TA \vee . \sim TB \vee TA) \rightarrow.$$
By Theorem 5.13, $\sim (\sim TA \vee . \sim TB \vee TA) \rightarrow$ is thought of as 
$F[TA^+, TA^-] \rightarrow$. (Two identical formulas are obviously equivalent to each other in our sense. And 
$TA$ is equialent to $TA$.) $F[TA~+, TA^-] \rightarrow$, moreover, has no implication sign as remarked above. 
$F[TA^+, TA^-] \rightarrow$ is, therefore, provable in $\bf PCN$ by Lemma 7.1. 

Case 2 $(\mbox{\bf CN} \vdash^0 A \supset B . \supset . (A \supset . B \supset C) \supset (A \supset C))$: The sequent 
obtained by the translation is of the form;
$$\sim (\sim TA \vee TB) \vee . \sim (\sim TA \vee . \sim (\sim TB \vee TC)) \vee \sim (\sim TA \vee TC) \rightarrow.$$ 
We, then, have the reduction:
$$\displaystyle{\frac{\enspace \sim (\sim (\sim TA \vee TB) \vee . \sim (\sim TA \vee . \sim (\sim TB \vee TC)) 
\vee \sim (\sim TA \vee TC)) \rightarrow \enspace \enspace }
{\enspace S_1 \enspace \enspace \enspace | \enspace \enspace 
\displaystyle{\enspace \enspace \frac{\enspace \enspace S_2 \enspace \enspace }
{\enspace \enspace \enspace S_3 \enspace \enspace \enspace | \enspace \enspace 
\displaystyle{\enspace \enspace \frac{\enspace \enspace S_4 \enspace \enspace }
{\enspace \enspace \enspace S_5 \enspace \enspace \enspace | \enspace \enspace 
\enspace \enspace S_6 \enspace \enspace , } \enspace (\vee^+ \rightarrow)}
}\enspace (\vee^+ \rightarrow)} 
}} \enspace (\vee^+ \rightarrow)$$
where 

\medskip

$S_1 = \sim (\sim \sim TA \vee . \sim (\sim TA \vee . \sim TB \vee TC) \vee (\sim TA \vee TC)) \rightarrow$, 

$S_2 = \sim (\sim TB \vee . \sim (TA \vee . \sim TB \vee TC) \vee (\sim TA \vee TC)) \rightarrow$, 

$S_3 = \sim (\sim TB \vee . \sim \sim TA \vee (\sim TA \vee TC)) \rightarrow$, 

$S_4 = \sim (\sim TB \vee . \sim (\sim TB \vee TC) \vee (\sim TA \vee TC)) \rightarrow$, 

$S_5 = \sim (\sim TB \vee . \sim \sim TB \vee (\sim TA \vee TC)) \rightarrow$, 

$S_6 = \sim (\sim TB \vee . \sim TC \vee (\sim TA \vee TC)) \rightarrow$. 

\medskip 

\noindent $S_1$, $S_3$, $S_5$ and $S_6$ are, respectively, thought of as:

\medskip 

$F[TA^+, TA^-] \rightarrow$, 

$F[TA^+, TA^-] \rightarrow$, 

$F_3[TB^+, TB^-] \rightarrow$, 

$F_4[TC^+, TC^-] \rightarrow$. 

\medskip 

\noindent By Theorem 5.13, $S_1$, $S_3$, $S_5$ and $S_6$ are all provable in $\bf SN$. 
They are also seen provable in $\bf PCN$ by taking their structures into consideration 
as in the case 1. The given sequent obtained by the translation is, thus, reduced to 
sequents provable in $\bf PCN$ in the above reduction, which are easily turned into complete tableaus. 
Thus, 
$$\bf PCN \it \vdash \enspace \sim (\sim (\sim TA \vee TB) \vee . \sim (\sim TA \vee . 
\sim TB \vee TC) \vee (\sim TA \vee T)) \rightarrow.$$
For brevity we will mention only the sequents obined by the translation and corresponding reductions 
in the cases 3--12 below. The detainls are similar to those of the cases 1 and 2.

Case 3 $(\mbox{\bf CN} \vdash^0 A \supset . B \supset (A \wedge B))$: Upon translation we have: 
$$\sim TA \vee . \sim TB \vee (TA \wedge TB) \rightarrow .$$
The reduction is, then, 
$$\displaystyle{\enspace \enspace \frac{\enspace \enspace \sim (\sim TA \vee . \sim TB \vee (TA \wedge TB)) \rightarrow \enspace \enspace }
{\enspace \enspace S_1 \enspace \enspace | 
\enspace \enspace S_2 \enspace \enspace }\enspace \enspace (\wedge^- \rightarrow)}$$
where 

$S_1 = \sim (\sim TA \vee . \sim TB \vee TA) \rightarrow$, 

$S_2 = \sim (\sim TA \vee . \sim TB \vee TB) \rightarrow$. 

\noindent $S_1$ and $S_2$ are, then, thought of as $F_1[TA^+, TA^-] \rightarrow$ and $F_1[TB^+, TB^-] \rightarrow$, respectively. 

Case 4 $(\mbox{\bf CN} \vdash^0 A \wedge B. \supset A)$: The sequent to be proved is in $\bf  PCN$ is 
$\sim (\sim (TA \wedge TB) \vee TB) \rightarrow$. The sequent is of the structure $F[TA^+, TA^-] \rightarrow$, and provable in $\bf PCN$. 

Case 5 $(\mbox{\bf CN} \vdash^0 A \wedge B. \supset A)$: The sequent obtained upon translation is 
$$\sim (\sim (TA \wedge TB) \vee TB) \rightarrow \enspace ,$$
which is of the form $F[TB^+, TB^-] \rightarrow$. 

Case 6 $(\mbox{\bf CN} \vdash^0 B \supset . A \vee B)$: As a result of translation we have $\sim (\sim TA \vee . TA \vee TB) \rightarrow$. 
This is thought of $F[TB^+, TB^-] \rightarrow$. 

Case 7 $(\mbox{\bf CN} \vdash^0 B \supset . A \vee B)$: This case is similar to the case 6. 

Case 8 $(\mbox{\bf CN} \vdash^0 A \supset C . \supset . (B \supset C . \supset (A \wedge B . \supset C))$: It is transformed into: 
$$\sim (\sim (\sim TA \vee TC) \vee . \sim (\sim TB \vee TC) \vee (\sim (TA \vee TB) \vee TC)) \rightarrow.$$
The reduction proceeds in the following way: 
$$\displaystyle{\frac{\enspace \sim (\sim (\sim TA \vee TC) \vee . \sim (\sim TB \vee TC) \vee (\sim (TA \vee TB) \vee TC)) \rightarrow \enspace }
{
\displaystyle{\frac{\enspace \enspace S_1 \enspace \enspace }
{\enspace \enspace S_3 \enspace \enspace | 
\enspace \enspace S_4 \enspace \enspace } \enspace \enspace (\vee^+ \rightarrow) } \enspace \enspace 
 | \enspace \enspace 
\displaystyle{\frac{\enspace \enspace S_2 \enspace \enspace }
{\enspace \enspace S_5 \enspace \enspace | \enspace \enspace 
 S_6 \enspace \enspace  } 
}\enspace \enspace (\vee^+ \rightarrow)}\enspace \enspace (\vee^+ \rightarrow)}$$
where 

\medskip

$S_1 = \sim (\sim (\sim TA \vee TC) \vee . \sim (\sim TB \vee TC) \vee (\sim TA \vee TC)) \rightarrow$, 

$S_2 = \sim (\sim (\sim TA \vee TC) \vee . \sim (\sim TB \vee TC) \vee (\sim TB \vee TC)) \rightarrow$, 

$S_3 = \sim (\sim \sim TA \vee . \sim (\sim TB \vee TC) \vee (\sim TA \vee TC)) \rightarrow$, 

$S_4 = \sim (\sim TC \vee . \sim (\sim TB \vee TC) \vee (\sim TA \vee TC)) \rightarrow$, 

$S_5 = \sim (\sim (\sim TA \vee TC) \vee. \sim \sim TB \vee (\sim TB \vee TC)) \rightarrow$, 

$S_6 = \sim (\sim (\sim TA \vee TC) \vee. \sim TC \vee (\sim TB \vee TC)) \rightarrow$.

\medskip

\noindent $S_3$, $S_4$,  $S_5$ and $S_6$ are, respectively, thought of as:

\medskip 

$F_1[TA^+, TA^-] \rightarrow$, 

$F_2[TC^+, TC^-] \rightarrow$, 

$F[TB^+, TB^-] \rightarrow$, 

$F[TC^+, TC^-] \rightarrow$. 

\medskip 

Case 9 $(\mbox{\bf CN} \vdash^0 A \supset B . \supset . (A \supset \sim B) \supset \sim A)$: The formula is transformed into:
$$\sim (\sim TA \vee TB) \vee . \sim (\sim TA \vee \sim TB) \vee \sim TA.$$
This is subject to the following series of reductions:

$$\displaystyle{\enspace \enspace \frac{\enspace \enspace \sim 
(\sim (\sim TA \vee TB) \vee . \sim (\sim TA \vee \sim TB) \vee \sim TA) \rightarrow \enspace \enspace }
{\enspace \enspace \enspace S_1 \enspace \enspace \enspace | \enspace \enspace 
\displaystyle{\enspace \enspace \frac{\enspace \enspace S_2 \enspace \enspace }
{\enspace \enspace \enspace S_3 \enspace \enspace \enspace | \enspace \enspace 
\enspace \enspace S_4 \enspace \enspace , } \enspace (\vee^+ \rightarrow)}
}\enspace (\vee^+ \rightarrow)} 
$$
where

\medskip

$S_1 = \sim (\sim \sim TA \vee . \sim (\sim TA \vee \sim TB) \vee \sim TA) \rightarrow$, 

$S_2 = \sim (\sim TAB \vee . \sim (\sim TA \vee \sim TB) \vee \sim TA) \rightarrow$, 

$S_3 = \sim (\sim TB \vee . \sim \sim TA \vee \sim TA) \rightarrow$, 

$S_4 = \sim (\sim TB \vee . \sim \sim TB \vee \sim TA) \rightarrow$.

\medskip

\noindent Here, $S_1$, $S_3$ and $S_4$ are, respectively, of the forms:

\medskip 

$F_1[TA^+, TA^-] \rightarrow$, 

$F_2[TA^+, TA^-] \rightarrow$, 

$F_3[TB^+, TB^-] \rightarrow$.

\medskip 

Case 10 $(\mbox{\bf CN} \vdash^0 \sim \sim A \supset A)$: It is transformed into the sequent 
$\sim (\sim \sim TA \vee TA) \rightarrow$, which is of the form $F_1[TA^+, TA^-] \rightarrow$. 

Case 11 $(\mbox{\bf CN} \vdash^0 \forall xA(x) \supset A(t))$: The translation of the formula results in 
the sequent:
$$\sim (\sim \forall xTA(x) \vee TA(t)) \rightarrow.$$
The sequent is subject to the following reduction: 
$$\displaystyle{\frac{\sim (\sim \forall xTA(x) \vee TA(t)) \rightarrow}
{\sim (\sim \forall xTA(x) \vee TA(t)) \wedge TA(t)\rightarrow \enspace .}} \enspace  
(\forall^+ \rightarrow)$$
The sequent obtained by the reduction is thought of as $F[TA(t)^+, TA(t)^-] \rightarrow$.  

Case 12 $(\mbox{\bf CN} \vdash^0 A(t) \supset \exists xA(x))$: The sequent to be proved is of the form:
the sequent:
$$\sim (\sim TA(t) \vee \exists xTA(x))) \rightarrow.$$
By applying $(\exists^- \rightarrow)$, the above sequent is reduced to 
$$\sim (\sim TA(t) \vee \exists xTA(x)) \wedge \sim TA(t) \rightarrow.$$
which is regarded as $F[TA(t)^+, TA(t)^-] \rightarrow$. 

Case 13 $(\mbox{\bf CN} \vdash^0 a' = b' \vee a = b)$: The sequent to be proved is for the form: 
$$\sim (\sim a' = b' \supset a = b) \rightarrow . \enspace  \enspace (*)$$
If we substitute $0^{(k)}$ and $0^{(l)}$ for $a$ and $b$, respectively in (*), then the sequent: 
$$\sim (\sim (0^{(k)})' = (0^{(l)})' \vee 0^{(k)} = 0^{(l)}) \rightarrow , \enspace  \enspace (**)$$
is obtained. Assume $k = l$. (++) is, then, the axiom 1, since $0^{(k)} = 0^{(l)}$ is an antecedentnegative part of 
(**), and a true constant prime formula. Assume $k \neq l$. (**) is, then, the axiom 2, since $ (0^{(k)})' = (0^{(l)})'$, 
i.e., $0^{(k + 1)} = 0^{(l + 1)}$ is an antecedent positive part of (**) and a false constant prime formula. From the 
above, (*) is the axiom 8 since $k$ and $l$ are thought of as arbitrary natural numbers. 

Case 14 $(\mbox{\bf CN} \vdash^0 \sim  a' = 0)$:The sequent in question has the form: 
$$\sim \sim  a' = 0 \rightarrow . \enspace  \enspace (*)$$
If $0^{(k)}$ is substitued for $a$ in (*), then the sequent 
$$\sim \sim (0^{(k)})' = 0 \rightarrow , \enspace  \enspace (**)$$
is obtained. (*) is, the axiom 8 since (**) is the axiom 2 for every natural number $k$. More specifically, 
$(0^{(k)})' = 0$ is a false constant prime formula and an antecedent positive part of (**). 

Case 15 $(\mbox{\bf CN} \vdash^0 a = b \supset . a = c \supset b = c)$: It is transformed into:
$$\sim (\sim a = b \vee . \sim a = c \vee b = c) \rightarrow , \enspace  \enspace (*)$$
Let $A(x)$ be $x = c$. By Theorem 5.9.(1),(*) is provable in $\bf SN$ since the sequent is regarded as 
$F[a = b ^+, A(a)^+, A(b)^-] \rightarrow$. 

Case 16 $(\mbox{\bf CN} \vdash^0 a = b \supset a' = b')$: The sequent transformed is, then, the axiom 8 since we can 
almost similarly prove it as in the case 14 of the vasis case. 

Case 17 $(\mbox{\bf CN} \vdash^0 a + 0 = a)$: The sequent to be proved is obviously $\sim a+ 0 = a \rightarrow$. 
For every natural number $k$, $\sim 0^{(k)} + 0 = 0^{(k)} \rightarrow$ is the axiom 1 since $ 0^{(k)} + 0 = 0^{(k)}$ 
is antecedent negative part of the sequent and a true constant prime formula, 

Case 18 $(\mbox{\bf CN} \vdash^0 a + b' = (a + b)')$: The sequent to be proved is 
$$\sim  a + b' = (a + b)' \rightarrow . \enspace \enspace (*)$$
This satisfies the condition of the axiom 8 since
$F[0^{(k)} + (0^{(l)})' = (0^{(k)} + 0^{(l)})'^-] \rightarrow$, i.e., a substitution instance of (*) is the axiom 1. Like the 
treatment of the case 17, $0^{(k)} + (0^{(l)})' = (0^{(k)} + 0^{(l)})'$ is a true constant prime formula, since 
$0^{(k)} + (0^{(l)})'$ and $(0^{(k)} + 0^{(l)})'$ have the same value $k + l + 1$. 

Case 19 $(\mbox{\bf CN} \vdash^0 a \cdot 0 = 0)$: The translation leads to the sequent of the form 
$\sim a \cdot 0  = 0 \rightarrow$. The sequent is, then, the axiom 8 since $\sim 0^{(k)} \cdot 0 = 0 \rightarrow$ 
is the axiom 1 for every natural number $k$. 

Case 20 $(\mbox{\bf CN} \vdash^0 a \cdot b' = a \cdot b + a)$: The formula is transformed into 
$\sim a \cdot b' = a \cdot b + a \rightarrow$, which is the axiom8. In fact, given, particular natural number 
$k$ and $l$, we can decide that $0^{(k)} \cdot (0^{(l)})'$ and  $0^{(k)} \cdot 0^{(l)} + 0^{(k)}$ have the same 
value $k \cdot (l + 1)$. 

By the basis, we obtain the proposition
$$\mbox{\bf CN} \vdash A \enspace \enspace \Rightarrow \enspace \enspace \mbox{\bf SN} \vdash \sim TA \rightarrow .$$
Further, we have:
$$\mbox{\bf SN} \vdash \sim TA \rightarrow \enspace \enspace \Rightarrow \enspace \enspace \mbox{\bf PCN} \vdash \sim TA \rightarrow ,$$
since $TA$ contains no implication sign. This leads to the proposition:
$$\mbox{\bf CN} \vdash A \enspace \enspace \Rightarrow \enspace \enspace \mbox{\bf PCN} \vdash \sim TA \rightarrow .$$ 

\it Induction steps \rm $(n \geq 0)$ 

Case 1: Assume $A$ is obtained by the following inference (namely, modus ponens):
$$\displaystyle{ \enspace \enspace \frac{\mbox{\bf CN} \vdash^{n_1} B \enspace \enspace 
\mbox{\bf CN} \vdash^{n_2} B \supset A  \enspace \enspace }
{\mbox{\bf CN} \vdash^n A \enspace ,}}$$
where $n = max(n_1, n_2) +1$. By H.I. ($n_1 < n$ and $n_2 < n$), we obtain: 
$$\mbox{\bf PCN} \vdash \enspace \sim TB \rightarrow \enspace \enspace (*)$$ 
and 
$$\mbox{\bf PCN} \vdash \enspace \sim T(B \supset A) \rightarrow$$ 
i.e.,
$$\mbox{\bf PCN} \vdash \enspace \sim (\sim TB \vee TA) \rightarrow \enspace . \enspace (**)$$ 
Here, use is made of an easy lemma to the effect: 
$$\mbox{\bf PCN} \vdash A \rightarrow \enspace \enspace \Rightarrow 
\enspace \enspace \mbox{\bf FN} \vdash A \rightarrow \enspace ,$$  
by means of which the cut elimination theorem proved earlier turns into its counterpart in $\bf PCN$. By the 
$\bf PCN$-cut elimination theorem (Theorem 6.4), there obtains the looked-for sequent $\sim TA \rightarrow$ from (*) and (**). 

Case 2: Let $A$ be $C \supset \forall xA_1(x)$. Assume that the formula is inferred by way of arule of classical predicate logic in 
the following way:
$$\displaystyle{ \enspace \enspace \frac{\enspace  \enspace \mbox{\bf CN} \vdash^{n - 1} C \supset A_1(x) \enspace  \enspace }
{\enspace  \enspace \mbox{\bf CN} \vdash^n  C \supset \forall xA_1(x) \enspace .}}$$
By H.I., we obtain $ \mbox{\bf PCN} \vdash^n \sim T(C \supset \forall xA_1(x)) \rightarrow$, to which 
$\sim (\sim TC \vee \forall xTA_1(x)) \rightarrow$ is reduced by applying $(\forall^- \rightarrow)$ under the restriction on variables. 
$\sim (\sim TC \vee \forall xTA_1(x)) \rightarrow$ is obviously the sequent obtain by the translation of 
$C \supset \forall xA_1(x)$.

Case 3: Let $A$ be $\exists xA_1(x) \supset C$. Suppose that the formulais inferred by way of another rule of predicate logic 
as follows: 
$$\displaystyle{ \enspace \enspace \frac{\enspace  \enspace \mbox{\bf CN} \vdash^{n - 1} A_1(x) \supset C \enspace  \enspace }
{\enspace  \enspace \mbox{\bf CN} \vdash^n \exists xA_1(x) \supset C \enspace .}}$$
By H.I., we obtain $ \mbox{\bf PCN} \vdash^n \sim T(A_1(x) \supset C) \rightarrow$, i.e., 
$$\mbox{\bf PCN} \vdash^n \sim (\sim TA_1(x) \vee TC) \rightarrow \enspace . \enspace \enspace (*)$$
The looked-for sequent $\sim (\sim \exists xTA_1(x) \vee TC) \rightarrow$ is reduced to (*) by $(\exists^+ \rightarrow)$ under the 
restriction on variables. 

\it Proof of $\Leftarrow$ in the embedding theorem\rm . Assume $\mbox{\bf SN} \vdash \enspace \sim TA \rightarrow$. By Theorem 7.2, 
$\mbox{\bf SNH} \vdash \enspace \sim TA \supset \enspace \sim \sim TA$. Then, 
$\mbox{\bf CN} \vdash \enspace \sim TA \supset \enspace \sim \sim TA$ since $\bf SNH$ is 
a subsystem of $\bf CN$. $TA$ is, thus, provable in $\bf CN$ since  
$\mbox{\bf CN} \vdash \enspace \sim TA \supset \sim \sim TA . \equiv TA$. 
(Here, $A \equiv B$ stands for $A \supset B . \wedge . B \supset A$.) $A$ is, therefore, provable in $\bf CN$ since 
$\mbox{\bf CN} \vdash A \equiv TA$. $\Box$








 





\section{Consistency of $\bf CN$}

\begin{theorem} \rm (Consistency theorem of $\bf CN$) \it 

\medskip 

\enspace \enspace $\bf CN$ is  consistent. 

\end{theorem}

\it Proof. \rm Assume, if possible, $\bf CN$ were inconsistent. Every formula of $\bf CN$ would, then, be provable in $\bf CN$. 
$\mbox{\bf CN} \vdash 0 = 1$, thus. (We abbreviate $0'$ as $1$.) $\mbox{\bf PCN}  \vdash \enspace \sim 0 = 1 \rightarrow \enspace $ is 
forthcoming right away by the embedding theorem (Theorem 6.3). It is not provable in $\bf PCN$, however. In fact, 
$\sim 0 = 1 \rightarrow$ is not any axiom of $\bf PCN$ and is by no means reduced by applying any reduction rule of 
$\bf PCN$. 

Here, we, are presenting another consistencey proof of theorem. Assume $\bf CN$, if possible, be inconsistent. Any formula 
would, then, provable in $\bf CN$. We, thus, obtain $\mbox{\bf CN} \vdash  0 = 1$ and $\mbox{\bf CN} \vdash \enspace \sim 0 = 1$. 
Then, $\sim 0 = 1 \rightarrow \enspace  $ and $\sim \sim 0 = 1 \rightarrow \enspace$ would be both provable in 
$\mbox{\bf PCN}$ by the embedding theorem. It follows that $\mbox{\bf PCN} \vdash \rightarrow 0 = 1$ and 
$\mbox{\bf PCN} \vdash \enspace \rightarrow \enspace \sim 0 = 1$ by Theorem 5.12. This is 
impossible since $\bf SN$ is consistent (Theorem 6.2).  Another standard argument for the consistency is the following. By the 
$\bf SN$-cut elimination theorem, $\enspace \rightarrow \enspace$ is obtained from $\mbox{\bf SN} \vdash \rightarrow 0 = 1$ 
and $\mbox{\bf SN} \vdash \enspace \rightarrow \enspace \sim 0 = 1$. But $\enspace \rightarrow \enspace $ is by no means 
provable in $\bf SN$.  $\Box$

\section{The proof of $\bf SN$-cut elimination theorem and its restricted version}

We first recall our cut elimination theorem. 

\medskip 

\noindent \bf Theorem 6.4 \rm ($\bf SN$-cut elimination theorem) \it 
$$(\mbox{\bf SN}  \vdash \Gamma \rightarrow G[A_\pm]  \mbox{ and  }  \mbox{\bf SN} \vdash F[A^\pm] \rightarrow \Delta) 
\enspace \Rightarrow \enspace \mbox{\bf SN} \vdash \Gamma \wedge F[\; \;^\pm] \rightarrow G[\; \;_\pm] \vee \Delta \enspace.$$

\it Proof\rm . A syntactic proof may be taken as an adaptation of Kanai [1984] \cite{kanai} to our system as that of $\bf PCN$-cut 
elimination theorem.  A semantical proof is also possible. $\Box$

In this section, as a novelty, we shall here give a simple proof of a restricted version of $\bf SN$-cut elimination theorem as an 
application of the disjunction property, using  $\bf PCN$-cut elimination theorem. 

\begin{theorem} \rm (A restricted $\bf SN$-cut elimination theorem) \it We have: 
$$(1) \enspace (\mbox{\bf SN}  \vdash \Gamma \rightarrow G[A_+]  \mbox{ and  }  \mbox{\bf SN} \vdash F[A^+] \rightarrow \Delta) 
\enspace \Rightarrow \enspace \mbox{\bf SN} \vdash \Gamma \wedge F[\; \;^+] \rightarrow G[\; \;_+] \vee \Delta \enspace ,$$
where 
$$\mbox{\bf SN}  \vdash \Gamma \wedge \sim A \rightarrow \enspace \Rightarrow \enspace \mbox{\bf PCN}  \vdash \Gamma \wedge \sim A \rightarrow$$
and 
$$\mbox{\bf SN} \vdash F[A^+] \wedge \sim \Delta \rightarrow \enspace \Rightarrow \enspace \mbox{\bf PCN} \vdash F[A^+] \wedge \sim \Delta \rightarrow$$
hold.
$$(2) \enspace (\mbox{\bf SN}  \vdash \Gamma \rightarrow G[A_-]  \mbox{ and  }  \mbox{\bf SN} \vdash F[A^-] \rightarrow \Delta) 
\enspace \Rightarrow \enspace \mbox{\bf SN} \vdash \Gamma \wedge F[\; \;^-] \rightarrow G[\; \;_-] \vee \Delta \enspace ,$$
where 
$$\mbox{\bf SN}  \vdash \Gamma \wedge A \rightarrow \enspace \Rightarrow \enspace \mbox{\bf PCN}  \vdash \Gamma \wedge A \rightarrow$$
and 
$$\mbox{\bf SN} \vdash F[A^-] \wedge \sim \Delta \rightarrow \enspace \Rightarrow \enspace \mbox{\bf PCN} \vdash F[A^-] \wedge \sim \Delta \rightarrow$$
hold.
\end{theorem}

\it Proof\rm . First we shall prove (1). 
$$(\mbox{\bf SN}  \vdash \Gamma \rightarrow G[A_+]  \mbox{ and  }  \mbox{\bf SN} \vdash F[A^+] \rightarrow \Delta) 
\enspace \Rightarrow \enspace \mbox{\bf SN} \vdash \Gamma \wedge F[\; \;^+] \rightarrow G[\; \;_+] \vee \Delta \enspace.$$
Suppose $\mbox{\bf SN}  \vdash \Gamma \rightarrow G[A_+]  \mbox{ and  }  \mbox{\bf SN} \vdash F[A^+] \rightarrow \Delta$.  
By the disjunction property (Theorem 5.7), we have $\mbox{\bf SN}  \vdash \Gamma \rightarrow G[\; \;_+]$ or 
$\mbox{\bf SN}  \vdash \Gamma \rightarrow A$ Assume $\mbox{\bf SN}  \vdash \Gamma \rightarrow G[\; \;_+]$. Then we 
immediately obtain $\mbox{\bf SN} \vdash \Gamma \wedge F[\; \;^+] \rightarrow G[\; \;_+] \vee \Delta$ by the thinning theorem. 
Assume $\mbox{\bf SN}  \vdash \Gamma \rightarrow A$. By Theorem 5.10, we get 
$\mbox{\bf SN}  \vdash \Gamma \wedge \sim A \rightarrow$ and $\mbox{\bf SN} \vdash F[A^+] \wedge \sim \Delta \rightarrow$.
From the assumption of the theorem, we obtain 
$$\mbox{\bf PCN}  \vdash \Gamma \wedge  \sim A \rightarrow , \enspace \enspace (*1)$$
$$\mbox{\bf PCN} \vdash F[A^+] \wedge \sim \Delta \rightarrow . \enspace \enspace (*2)$$
Apply $\bf PCN$-cut elimination theorem to $(*1)$ and $(*2)$. Then we obtain 
$$\mbox{\bf PCN} \vdash \Gamma \wedge F[\; \;^+] \wedge \sim \Delta \rightarrow . \enspace \enspace (*3)$$
From $(*3)$ and Theorem 5.11, we have $\mbox{\bf SN} \vdash \Gamma \wedge F[\; \;^+] \rightarrow \Delta$. By the thinning theorem, 
the disired sequent 
$$\mbox{\bf SN} \vdash \Gamma \wedge F[\; \;^+] \rightarrow G[\; \;_+] \vee \Delta$$
holds. For the proof of (2), we can take similar arguments. $\Box$

The restriction required by Theorem 9.1 is very strong. But some applications of Theorem 9.1 seem to be possible in computer science, 
when we pursue to constructive nature appeared in the subject.

\section{On future studies}

There may be a plan to extend the resuls of this paper to the full system with complete induction, a system with bar induction and formal 
analysis, etc. For those studies, the following papers will be useful: Bezem [1985, 1989] \cite{bezem1, bezem2}, Pohlers [2009] \cite{pohlers}, 
Ferreira [2015] \cite{fer}, Kahle and Rathjen [2015] \cite{kr}, Siders [2015] \cite{siders}, Spector [1962]  \cite{spector}, Tait [2015] \cite{tait}, etc. 

It would also be a way to build the strong negation version of Sch\"utte's book, Proof Theory (1977 version) \cite{schutte1977}. I shall have 
a plan to do so for the future Part II of this paper. 


\section*{Ackowledgements}

I would like to thank the late Professor Emeritus Arata Ishimoto for introducing me strong negation and Novikov's proof theory (regularity), and his encouragements 
to me during our long joint research. This paper is an extension and a corrected refinement of the appendix of my Master thesis (Tokyo University of Science, 1984,) 
under his supervision. 

To write this paper, I was very much inspired by the papers and the book, Siders [2015]  \cite{siders}, Ferreira [2015] \cite{fer} and 
Odintsov [2008] \cite{odin}. So, I would like to thank Prof. Siders, Prof. Ferreira and Prof. Odintsov. I also appreciate some reviews of 
zbMATH by Prof. M. Yasuhara and Prof. M. J. Beeson. 



\noindent Meiji Pharmaceutical University

\noindent Department of Medical Molecular Informatics

\noindent Tokyo, Japan 

\bigskip

\noindent Hosei University

\noindent Graduate School of Science and Engineering

\noindent Tokyo, Japan 

\bigskip

\noindent Hosei University

\noindent Faculty of Science and Engineering

\noindent Department of Applied Informatics

\noindent Tokyo, Japan

\medskip

\noindent ta-inoue@my-pharm.ac.jp

\noindent takao.inoue.22@hosei.ac.jp

\noindent takaoapple@gmail.com


\begin{thebibliography}{99}

\bibitem{aczel} \sc P. H.  G.~Aczel, H. Simmons \rm and \sc S. S. Wainer\rm , (eds.) , {\bf Proof Theory}\rm , Cambridge U. P.,1992.

\bibitem{ackermann1924} \sc W.~Ackermann\rm , {Begr\"undung des ‘Tertium non datur’mittels der Hilbertschen Theorie der Widerspruchsfreiheit}, \bf Mathematische 
Annalen\rm , Vol. 93 (1924), pp. 1--36.

\bibitem{ackermann1940} \sc W.~Ackermann\rm , {Zur Widerspruchsfreiheit der Zahlentheorie}, \bf Mathematische Annalen\rm , Vol. 117 (1940), pp. 162--194.

\bibitem{akama} \sc S. Akama\rm , {Constructive predicate logic with strong negation and model theory}, \bf Nortre Dame Journal of Formal Logic\rm , Vol.29 (1988),  pp. 18--27.

\bibitem{an} \sc A. Almukdad \rm and \sc D. Nelson\rm , {Constructible falsity and inexact predicates}, \bf The Journal of Symbolic Logic\rm , Vol. 49 (1984), pp. 231--233.

\bibitem{arai} \sc T. Arai\rm , {\bf Mathematical Logic}\rm , (in Japanese), Iwanami Shoten, Tokyo, 2011. 

\bibitem{bellotti} \sc L. Bellotti\rm , {Novikov's cut elimination}, \bf Logique et Analyse\rm , Vol.61, No. 242 (2018), pp. 183-199. 

\bibitem{bernays} \sc P. Bernays\rm , {On the original Gentzen consistency proof for number theorey}, In \sc A. Kino, J. Myhill \rm and \sc R. E. Vesley\rm (eds.), 
\bf Intuitionism and Proof Theory\rm , pp. 409--417, North-Holland, Amsterdam, 1970. 

\bibitem{bezem1} \sc M. Bezem\rm , Strongly majorizable functionals of finite type: a model for bar recursion containing
discontinuous functionals, \bf The Journal of Symbolic Logic\rm , Vol.50 (1985), pp.652--660.

\bibitem{bezem2} \sc M. Bazem\rm , Compact and majorizable functionals of finite type, \bf The Journal of Symbolic Logic\rm , Vol.54 (1989), pp.271--280.

\bibitem{buch} \sc W. Buchholz, S. Feferman, W. Pohlers \rm and \sc W. Sieg, \bf Iterated Inductive Definitions and Subsystems of Analysis: Recent Proof-Theoretical Studies\rm (Lecture Notes in Mathematics Vol. 897),  Springer, 1981. 

\bibitem{busshandbook} \sc S. L. Buss\rm , (ed.), {\bf Handbook of Proof Theory}\rm , Elsevier, Amsterdam, 1998.

\bibitem{bussIg} \sc S. L. Buss \rm and \sc A. Ignjatovi'c\rm , Unprovability of Consistency Statements in Fragments of Bounded Arithmetic, \bf Annals of Pure and Applied Logic\rm , Vol.74 (1995), pp. 221-244.

\bibitem{dekker} \sc J. C. E. Dekker \rm (ed.), \bf Recursive function theory\rm , Proceedings of Symposia in Pure Mathematics, Vol. 5, American Mathematical Society, Providence, Rhode Island, 1962.

\bibitem{feferman} \sc S. Feferman\rm , {How we got from there to here}, In \cite{buch}, pp. 1--15.

\bibitem{fs1} \sc S. Feferman  \rm and \sc W. Sieg, \rm Iterated inductive definitions and subsystems of analysis. In \cite{buch}, pp. 16--77. 

\bibitem{fs2} \sc S. Feferman  \rm and \sc W. Sieg, \rm Proof theoretic equivalences between classical and constructive theories for analysis. In \cite{buch}, pp. 78--142. 

\bibitem{fer} \sc F. Ferreira\rm , Spector's Proof of the Consistency of Analysis, In \cite{kr}, pp. 279--300. 

\bibitem{ff} \sc R. C. Flagg \rm and \sc H. Friedman\rm , {Epistemic and intuitionistic formal systems}, \bf Annals of Pure and Applied Logic\rm , Vol. 32 (1986), pp. 53--60.

\bibitem{fitting} \sc M. Fitting\rm , {\bf Proof Methods for Modal and Intuitionistic Logics}, D. Reidel, Dordrecht, 1983.

\bibitem{friedf} \sc R. C. Flagg \rm and \sc H. Friedman\rm , {A framework for measuring the complexity of mathematical concepts}, \bf Advances in Pure and Applied Mathematics\rm , Vol. 11 (1991), pp. 1--34.

\bibitem{gentzen1934a} \sc G. Gentzen\rm , {Untersuchungen \"uber das Logische Schliessen I}, \bf Math. Zeitschrift\rm , Vol. 39 (1934), pp. 176--210.

\bibitem{gentzen1934b} \sc G. Gentzen\rm , {Untersuchungen \"uber das Logische Schliessen II}, \bf Math. Zeitschrift\rm , Vol. 39 (1934), pp. 405--431. 

\bibitem{gentzen1936} \sc G. Gentzen\rm , {Wiederspruchsfreiheit der reinen Zahlentheorie}, \bf Math. Ann.\rm , Vol. 112 (1936), pp. 493--565.

\bibitem{gedel} \sc K. G\"{o}del\rm , {Zur intuitionistischen Arithmetik und Zahlentheorie}, \bf Ergebnisse eines mathematischen Kolloquiums\rm , Vol. 4 (1933), pp. 34--38. 

\bibitem{herbrand} \sc J. Herbrand\rm , {Sur la non-contradiction de l'arithm\'etique}, \bf Journal f\"ur reine und angewandte Mathematik\rm , Vol. 166 (1931), 
pp. 1--8. English translation: ‘On the consistency of arithmetic’,  In \sc J. van Heijenoort \rm(ed.) , \bf From Frege to G\"odel\rm , 
Harvard U. P. , Cambridge, 1967, pp. 620--628. 

\bibitem{hb1970} \sc D. Hilbert \rm and \sc P. Bernays, {\bf Grundlagen der Mathematik, Bd. I (1934), Bd. II (1939)}\rm , Springer, Berlin. 
(Second edition. 1968 and 1970, respectively)

\bibitem{inouen} \sc T. Inou\'{e}\rm , {New type cut elimination theorem and its application to embedding theorem in S}. Unpublished note, 1984.

\bibitem{inoueac} \sc T. Inou\'{e}\rm , {\rm On the Accommodation of Prolog to Constructive Logic with Strong Negation S using its Embedding 
Theorem (in Japanese)},  In \cite{ishimotoprog}, pp. 41--93. 

\bibitem{inouefunda} \sc T. Inou\'{e}\rm , {\rm On the fundamental theorem of Logic programming: the reduction of Prolog to constructive Logic with 
strong negation}, in preparation. 

\bibitem{ishimotoe} \sc A. Ishimoto\rm , {Constructive propositional logic with strong negation and their completeness}. Unpublished note, 197?.

\bibitem{ishimotos} \sc A. Ishimoto\rm , {A Sch\"utte-type formulation of the intuitionistic 
functional calculus with strong negation}, \bf Bulletin of the Tokyo Institute of Technology\rm , Vol. 100 (1970), pp. 161--189. 

\bibitem{ishimotoprog} \sc A. Ishimoto\rm , (ed.), \bf The progress report III of Ishimoto's group\rm , Grant-in-Aid for Scientific Research No.57115012, 
Ministry of Education, Japan (1984),

\bibitem{ishimototaga} \sc A. Ishimoto\rm , (ed.), \bf The Logic of Natural Language and its Ontology \rm  (in Japanese), Taga Shuppan, Tokyo, 1990.

\bibitem{kr} \sc R. Kahle \rm and \sc M. Rathjen \rm (eds.), \bf Gentzen's Centenary, The Quest for Consistency\rm , Springer, 2015. 

\bibitem{kanai} \sc N. Kanai\rm , {The cut elimination theorem of constructive predicate logic with strong negation (in Japanese)}. 
In \cite{ishimotoprog}, pp. 94--105.

\bibitem{khlo} \sc I. N. Khlodovskii\rm , A new proof of the consistency of arithmetic, (in Russian)
\bf Uspekhi Mat. Nauk\rm , Vol.14 (1959), pp.105--140. (as \sc I. N. Holdovskii\rm ) English translation of the
preceding by Moshe Machover. American Mathematical Society translations, ser. 2, Vol. 23 
(\bf Nine Papers on Logic and Quantum Electrodynamics\rm ) (1963), pp. 191-230.

\bibitem{kleene0} \sc S. C. Kleene\rm , {Permutability of inferences in Gentzen's calculi LK and LJ}, in \sc S. C. Kleene\rm , \bf Two Papers on the Predicate 
Calculus\rm , \bf Memoirs of the American Mathematical Society\rm , No. 10 (1952), The American Mathematical Society, Providence, Rhode Island, pp. 1--26.

\bibitem{kleene1} \sc S. C. Kleene\rm , {\bf Introduction to Metamathematics}, North-Holland, Amsterdam, 1952. 

\bibitem{kleene2} \sc S. C. Kleene\rm , {\bf Mathematical Logic}, J. Wiley and Sons, New York, 1967.  (There is a Russian translation by G. E. Mints.)

\bibitem{kolmo} \sc A. N. Kolmogorov\rm , {On the principle of exculuded middle (in Russian)}. In \sc J. van Heijenoort \rm(ed.) , \bf From Frege to G\"odel\rm , pp. 417-437,
 Cambridge, Harvard U. P., 1967. 

\bibitem{maehara} \sc S. Maehara\rm , {\bf Mathematical Logic}, (in Japanese), Baihukan, Tokyo, 1973. 

\bibitem{markov} \sc A. A. Markov\rm , {A constructive logic (in Russian)}, \bf Uspehi Mathematiceskih\rm , Vol. 5 (1950), pp. 187--188. 

\bibitem{mi} \sc T. Matsuda \rm and \sc A. Ishimoto\rm , Relationship between logic computer language and constructive predicate logic with strong negation, (in Japanese),  
In \cite{ishimotoprog}, pp. 109--121. 

\bibitem{mura} \sc R. Murawski\rm , On proofs of the consistency of arithmetic, \bf Studies in Logic, Grammar and Rhetoric\rm , Vol. 4 (2001), pp.41--50.

\bibitem{nelson} \sc D. Nelson\rm , {Constructible falsity}, \bf The Journal of Symbolic Logic\rm , Vol. 14 (1949), pp. 16--26.

\bibitem{neumann} \sc J. von Neumann\rm , {Zur Hilbertschen Beweistheorie}, \bf Mathematische Zeitschrift\rm , Vol. 26 (1927), pp. 1--46.

\bibitem{monk} \sc D. Monk\rm , {\bf Mathematical Logic}, Springer-Verlag, New York, 1976.

\bibitem{novikov} \sc P. S. Novikov\rm , {\bf Elements of Mathematical Logic}, (in Russian), Moscow, 1959. (English translation: \bf Elements of 
Mathematical Logic \rm (translated by \sc L. F. Boron\rm ), Oliver \& Boyd, Edinburgh and London, 1964.) (German translation: \bf Grundz\"{u}ge 
der mathematischen Logik \rm (translated by \sc K. Rosenbaum\rm ), Friedr. Vieweg $+$ Sohn, Braunschweig, 1973. There is a Japanese translation 
by the late Prof. A. Ishimoto, Tokyo Tosho, Tokyo, 1965. 

\bibitem{odin} \sc S. P. Odintsov\rm , {\bf Constructive Negations and Paraconsistency}\rm , Springer, 2008.

\bibitem{onok} \sc K. Ono\rm , {Logische Untersuchungen \"{u}ber die Grundlagen der Mathematik}, \bf Journal of the Faculty of Science I, Imperial University of Tokyo\rm , Vol. 3 (1938), pp. 329--389. 
\bibitem{rasiowa} \sc H. N. Rasiowa\rm , Lattices and constructive logic with strong negation, \bf Fundamenta Mathematicae\rm , Vol.46 (1958), pp.61--80. 

\bibitem{pohlers} \sc W. Pohlers\rm , {\bf Proof Theory, The First Step into Impredicativity}\rm , Springer, Berlin, 2009. 

\bibitem{schutteB1960} \sc K. Sch\"{u}tte\rm , {\bf Beweistheorie}, Springer-Verlag, Berlin, 1960.

\bibitem{schutteI1968} \sc K. Sch\"{u}tte\rm , {Der Interpolationsatz der intuitionistischen Pr\"{a}dikatenlogik}, \bf Mathematishe Annalen\rm , Vol. 148 (1962), pp. 192--200.

\bibitem{schutteV1968} \sc K. Sch\"{u}tte\rm , {\bf Vollst\"{a}ndige Systeme Modaler und Intuitionistischer Logik}, Springer-Verlag, Berlin, 1968. 

\bibitem{schutte1977} \sc K. Sch\"{u}tte\rm , {\bf Proof Theory}, Springer-Verlag, Berlin, 1977. 

\bibitem{shimizu} \sc S. Shimizu\rm , A study on constructive logic with strong negation - its soundness and completeness, (in Japanese). In \cite{ishimototaga}, pp. 241--267. 

\bibitem{siders} \sc A. Siders\rm , A Direct Gentzen-Style Consistency Proof for Heyting Arithmetic. In \cite{kr}, pp.177--211.

\bibitem{smullyan} \sc R. M. Smullyan\rm, \bf First-Order Logic\rm , Springer, New-York, 1968. There is a Dover verion of this book with a short comment for literature. 

\bibitem{spector} \sc C. Spector\rm , Provably recursive functionals of analysis: a consistency proof of analysis by an extension of principles formulated in 
current intuitionistic mathematics. In \cite{dekker}, pp. pp. 1-27.

\bibitem{tait} \sc W. W. Tait\rm , Gentzen's original consistency proof and the Bar theorem, In \cite{kr}, pp. 213--228,

\bibitem{takano} \sc M. Takano\rm , {A formulation of the Fitch functional calculas as a Sequenzenkalk\"ul}, \bf Bulletin of the Tokyo Institue of Technology\rm , Vol. 100 (1970), pp. 143--160.

\bibitem{takeuti} \sc G. Takeuti\rm , {\bf Proof Theory\rm , 2nd}, North-Holland, Amsterdam, 1987. 

\bibitem{toledo} \sc S. Toledo\rm , {\bf Tableau Systems for First Order Number Theory and Certain Higher Order Theories}\rm, (Lecture Notes in Mathematics Vol. 447), 
Springer, Berlin, 1975.

\bibitem{ts} \sc A. S. Troelstra \rm and \sc H. Schwichtenberg\rm , {\bf Basic Proof Theory}, 2nd, Cambridge University Press, New York, 2000.

\bibitem{tvd} \sc A. S. Troelstra \rm and \sc D. van Dalen\rm , {\bf Constructivism in Mathematics, An Introduction, Vol. I, II}, North-Holland, Amsterdam, 1988. 

\bibitem{voro1} \sc N. N. Vorob'ev\rm , {A constructive propositional calculus with strong negation (in Russian)}, 
\bf Doklady Akademii Nauk SSSR\rm , Vol. 85 (1952), pp. 465--468.

\bibitem{voro2} \sc N. N. Vorob'ev\rm , {The problem of deducibility in constructive propositional
calculus with strong negation (in Russian)}, \bf Doklady Akademii Nauk SSSR\rm , Vol. 85 (1952), pp. 689--692. 

\bibitem{voro3} \sc N. N. Vorob'ev\rm , {Constructive propositional calculus with strong negation (in Russian)}, 
\bf Transactions of Steklov's Institute\rm , Vol. 72 (1964), pp. 195--227. 

\end{thebibliography}
\end{document}